\documentclass[12pt]{article}
\usepackage{amssymb,amsbsy,amsmath,amsfonts,amssymb}
\usepackage{latexsym,euscript,exscale}

\usepackage{times}

\newcommand{\pf}{

\smallskip

\noindent {\it Proof : }}

\newcommand {\N}{\mathbb N}

\newcommand {\R}{\mathbb R}

\newcommand {\Z}{\mathbb Z}

\newcommand {\C}{\mathbb C}
\newcommand {\HH}{\mathbb H}

\newcommand{\pff}{$\hfill \square$
\smallskip}

\newcommand{\norm}[1]{\ensuremath{\left\|#1\right\|}}

\newtheorem{prop}{Proposition}
\newtheorem{lemm}[prop]{Lemma}
\newtheorem{theo}[prop]{Theorem}
\newtheorem{ques}[prop]{Question}
\newtheorem{coro}[prop]{Corollary}

\newtheorem{defi}[prop]{Definition}

\newtheorem{exam}[prop]{Example}

\title{Countable groups of isometries on Banach spaces}

\author{Valentin Ferenczi and El\'oi  Medina Galego}

\date{ferenczi@ccr.jussieu.fr, eloi@ime.usp.br}

\begin{document}

\maketitle
 
\

\begin{abstract} A group $G$ is representable in a Banach space $X$ if $G$ is isomorphic to the group of isometries on $X$ in some equivalent norm. We prove that a countable group $G$ is representable in a separable real Banach space $X$ in several general cases, including when $G \simeq \{-1,1\} \times H$, $H$ finite and $\dim X \geq |H|$, or when $G$ contains a normal subgroup with two elements and $X$ is of the form
$c_0(Y)$ or $\ell_p(Y)$, $1 \leq p <+\infty$. This is a consequence of a result inspired by methods of S. Bellenot and stating that under rather general conditions on a separable real Banach space $X$ and a countable bounded group $G$ of isomorphisms on $X$ containing $-Id$, there exists an equivalent norm on $X$ for which $G$ is equal to the group of isometries on $X$.

We also extend methods of K. Jarosz to prove that any complex Banach space of dimension at least $2$ may be renormed to admit only trivial real isometries, and that any real Banach space which is a cartesian square may be renormed to admit only trivial and conjugation real isometries. It follows that every real Banach space of dimension at least $4$ and with a complex structure may be renormed to admit exactly two complex structures up to isometry, and that every real cartesian square may be renormed to admit a unique complex structure up to isometry.
\end{abstract}

\

\section{Introduction}

What groups $G$ may be seen as the group of isometries on a Banach space $X$? This general question may be formulated by the following definition given by K. Jarosz in 1988.

\begin{defi}(Jarosz \cite{J}) A group $G$ is representable in a Banach space $X$ if there exists an equivalent norm on $X$ for which the group of isometries on $X$ is isomorphic to $G$.
\end{defi}

In \cite{J}, Jarosz stated as an open question which groups were representable
in a given Banach space. The difference with the classical theory of
representation of groups on linear spaces is that here we require an
isomorphism with {\em the} group of isometries on a Banach space, and not just
some group of isometries or isomorphisms. Since $\{-Id,Id\}$ is always a normal subgroup of the group of isometries on a real Banach space, it follows that a group which is representable in a real Banach space must always contain a normal subgroup with two elements. Conversely, J. Stern \cite{St} proved that for any group $G$ which contains a normal subgroup with two elements, there exists a real Hilbert space 
$H$ such that $G$ is representable in $H$. Furthermore if $G$ is countable then $H$ may be chosen to be separable.

For an arbitrary Banach space $X$ it remains open which groups are representable in $X$. Jarosz proved that   $\{-1,1\}$ is representable in any real Banach space, and that the unit circle $C$ is representable in any complex space (the separable real case had been solved previously by S. Bellenot \cite{B}). He also proved that for any countable group $G$, $\{-1,1\} \times G$  is representable in $C([0,1])$, and that for any group $G$ there exists a complex space $X$ such that $C \times G$ is representable in $X$.

\

In a first section of this paper, we give a much more general answer to the question of representability by proving that:
\begin{itemize}
\item the group $\{-1,1\} \times G$ is representable in $X$ whenever $G$ is a finite group and $X$ a separable real space $X$ such that $\dim X \geq |G|$, Theorem \ref{fini},
\item the group $G$ is representable in $X$ whenever $G$ is a countable group admitting a normal subgroup with two elements and $X$ is a separable real Banach space with a symmetric decomposition either isomorphic to $c_0(Y)$ or to $\ell_p(Y)$ for some $Y$ and $1 \leq p <+\infty$, or with the Radon-Nikodym Property, Theorem \ref{classical},  
\item the group $\{-1,1\} \times G$ is representable in $X$
whenever $G$ is a countable group
and $X$ an infinite-dimensional separable real Banach space containing a complemented subspace with a symmetric basis, Theorem \ref{deno}.
\end{itemize}

These results are partial answers to a conjecture of Jarosz who asked whether $\{-1,1\} \times G$ is representable in $X$ for any group $G$ and any real space $X$ such that $\dim X \geq |G|$.

As an application of our results we obtain that a countable group $G$ is representable in $c_0$ (resp. $C([0,1])$, $\ell_p$ for $1 \leq p <+\infty$, $L_p$ for $1 \leq p<+\infty$) if and only if it contains a normal subgroup with two elements, Corollary \ref{Lp}.

\

Our method is to ask, given a group $G$ of linear isomorphisms on a real Banach space $X$, whether there exists an equivalent norm on $X$ for which $G$ is the group of isometries on $X$. Once the problem of representability is reduced to representing a given group as some group of isomorphisms on a given Banach space, it is  much simpler to address, and this leads to Theorem \ref{fini}, Theorem \ref{classical}, and Theorem \ref{deno}.  In other words, we explore in which respect the question of representability of groups in Banach spaces belongs to the renorming theory or rather may be reduced to the purely isomorphic theory.

If a group of isomorphisms is the group of isometries on a real (resp. complex) Banach space in some equivalent norm, then it must be bounded, contain $-Id$ (resp. $\lambda Id$ for all $\lambda \in C$), and be closed for the convergence of $T$ and $T^{-1}$ in the strong operator topology. Therefore the question is:

\begin{ques} Let $X$ be a real (resp. complex) Banach space and let $G$ be a group of isomorphisms on $X$ which is bounded, contains $-Id$ (resp. $\lambda Id$ for all $\lambda \in C$), and is closed for the convergence of $T$ and $T^{-1}$ in the strong operator topology.
Does there exist an equivalent norm on $X$ for which $G$ is the group of isometries on $X$?
\end{ques}

A positive answer was obtained by Y. Gordon and R. Loewy \cite{GL} when $X=\R^n$ and $G$ is finite; this answered a question by J. Lindenstrauss. In this paper, we extend the methods of Bellenot and use a renorming method of G. Lancien \cite{L} to considerably improve this result:

\begin{itemize} \item Let $X$ be a separable real Banach space with the Radon-Nikodym Property. Then for any countable bounded group $G$ of isomorphisms on $X$ which contains $-Id$ and is separated by some point with discrete orbit, there exists an equivalent norm on $X$ for which $G$ is equal to the group of isometries on $X$, Theorem \ref{avecRNP}.

\item Let $X$ be a separable real Banach space. Then for any finite group $G$ of isomorphisms on $X$ which contains $-Id$, there exists an equivalent norm on $X$ for which $G$ is equal to the group of isometries on $X$, Theorem \ref{renormingG}.
\end{itemize}

Therefore for separable real spaces and finite groups, the question of representability really does not belong to renorming theory. Also, note that a countable group of isomorphisms on $X$ which is equal to the group of isometries in some equivalent norm must always be discrete for the convergence of $T$ and $T^{-1}$ in the strong operator topology and admit a separating point, Lemma \ref{discrete}. It remains unknown however whether this implies the existence of a separating point with discrete orbit, that is, if the implication in Theorem \ref{avecRNP} is an equivalence for countable groups.

To conclude that section we deduce Theorem \ref{fini}, Theorem \ref{classical} and Theorem \ref{deno} essentially from Theorem
\ref{renormingG} and Theorem \ref{avecRNP}.
We also prove that Theorem \ref{fini} and Theorem \ref{deno} are optimal in the sense that there exists a real space in which representable finite groups are exactly those of the form $\{-1,1\} \times G$, Proposition \ref{HIreel}, and a real space containing a complemented subspace with a symmetric basis in which representable countable groups are exactly those of the form $\{-1,1\} \times G$, Proposition \ref{symmetricplusHI}. On the other hand we have the classical examples of $c_0$, $C([0,1])$, $\ell_p, 1 \leq p <+\infty$ and $L_p, 1 \leq p<+\infty$ for which Corollary \ref{Lp} states that representable countable groups are exactly those which admit a normal subgroup with two elements, and we also provide an intermediary example of a space in which the class of
representable finite groups is strictly contained in between the above two classes, Proposition \ref{HIcomplexe}.

\

In a second section of this paper, we use the renorming methods of Jarosz to study complex structures on  real Banach spaces up to isometry.
Our results are actually related to the representability of the circle group $C$ and of the group of isometries on $\C$ as the group of $\R$-linear isometries on a complex Banach space.

We recall a few facts about complex structures. Any complex Banach space is also a real Banach space, and conversely, the
linear structure on a real Banach space $X$
may be induced by a $\C$-linear structure; the corresponding complex Banach space is
said to be a {\em complex structure} on $X$ in the isometric sense.
  It is clear that any complex structure on $X$  is canonically associated
 to some  $\R$-linear map $I$ on $X$ such that $I^2=-Id$ and $\cos \theta Id+\sin \theta I$ is an isometry for all $\theta$, and which
 defines the multiplication by the imaginary number $i$. Conversely for any
 such map $I$, there exists an associated complex structure denoted
$X^I$.

The existing theory of complex structure, however, is up to isomorphism. In this case, complex structures correspond to real isomorphisms $I$ of square $-Id$, up to the renorming $\||.\||$ defined by $\||x\||=max_{\theta}\norm{\cos\theta x+\sin\theta Ix}$.
It is well-known that complex structures do not always exist up to isomorphism on a Banach space.
By  \cite{B}, \cite{K} there exists real spaces with at least two complex structures up to isomorphism, and the examples of \cite{B} and \cite{An} (which are separable) actually admit 
 a continuum of complex structures.
By \cite{F} for each $n \in \N^*$ there exists a space with exactly $n$
complex structures up to isomorphism. In \cite{F} and \cite{FG}  various examples of
spaces  different from the
classical example of $\ell_2$ are also shown to have a unique complex structure up
to isomorphism, including a HI example, a space with an
unconditional basis, and a $C(K)$ space defined by Plebanek.

It turns out actually that the classical spaces $c_0$, $C([0,1])$, $\ell_p, 1 \leq p
\leq +\infty$ and $L_p, 1 \leq p <+\infty$ also admit a unique complex
structure up to isomorphism. A nice and simple proof of this fact was given to us by
N.J. Kalton after a first version of our paper was posted and is included
here, Theorem \ref{kalton}.

\

The isometric theory of complex structures turns out to be totally different from the isomorphic theory. For a very general class of Banach spaces $X$, we show that quite various situations may be obtained concerning existence and uniqueness of complex structures up to isometry on $X$ by choosing different renormings on $X$. This may justify why the isometric theory of complex structures has not been investigated before, as it is unclear what other results one may want to obtain in that area.

 We first prove that $\ell_2$ has a unique complex structure up to isometry, Proposition \ref{l2}. On the other hand, since Jarosz \cite{J} showed that every real Banach space may be renormed to admit only trivial isometries (i.e. the only isometries are $Id$ and $-Id$), every real Banach space may be renormed not to admit complex structures in the isometric sense.

Extending the methods of Jarosz \cite{J} we prove 1) that any complex Banach space of dimension at least $2$ may be renormed to admit only trivial real isometries, and 2) that any complex Banach space which is real isomorphic to a cartesian square, and whose complex law is the canonical one associated to the decomposition as a square, may be renormed to admit only trivial and conjugation isometries, Corollary \ref{maintheorem}. It follows 1) that every real Banach space of dimension at least $4$ and with a complex structure may be renormed to admit exactly two complex structures up to isometry (the initial complex structure and its conjugate), and 2) that every real cartesian square may be renormed to admit as unique complex structure up to isometry the canonical complex structure associated to its decomposition as a square,
Theorem \ref{renormingcomplexstructure}.

\

In a last section we extend results of F. Rabiger and W.J. Ricker, \cite{RR}, by proving that any isometry on a real Banach space with the $\lambda Id+S$-property, such as Gowers-Maurey's space $X_{GM}$, is of the form $\pm Id+K$, $K$ compact, Proposition \ref{mnbvc}.

\

For classical results in Banach space theory, such as, for example, the definition of the Radon-Nikodym Property or of a symmetric basis, we refer to \cite{LT}; for renorming questions in Banach spaces, we refer to
\cite{DGZ}.

\section{Representation of countable groups on separable real Banach spaces}

\subsection{$G$-pimple norms on separable Banach spaces}
 
In this subsection we extend the construction of Bellenot in \cite{Bel} from $\{-Id,Id\}$ to countable groups of isometries.
So in the following, $X$ is real separable, $G$ is a countable group of isometries on $X$, and under certain conditions on $G$, we construct an equivalent norm on $X$ for which $G$ is the group of isometries on $X$.

Let us give an idea of our construction. Bellenot renorms $X$ with an LUR norm  and then defines, for $x_0$ in $X$ of norm $1$, a new unit ball
(the "pimple" ball) obtained by adding two small cones in $x_0$ and $-x_0$. Any isometry in the new norm must preserve the cones and therefore send $x_0$ to $\pm x_0$. Repeating this for a sequence $(x_n)_n$
with dense linear span, chosen carefully so that one can add the cones "independantly", and so that the sizes of the cones are "sufficiently" different, any isometry  sends $x_n$ to $\pm x_n$. Finally, if each $x_n$ was chosen "much closer" to $x_0$ than to $-x_0$, any isometry fixing $x_0$ must fix each $x_n$ and therefore any isometry is equal to $Id$ or $-Id$.

In our case one should obviously put cones of same size in each $gx_0, g \in G$, defining a "G-pimple ball"; therefore any isometry preserves the orbit $Gx_0$. Then one repeats a similar procedure as above, adding other cones in $gx_n, g \in G$ for a sufficiently dense sequence $(x_n)_n$, so that any isometry preserves $Gx_n$ for all $n$. These $x_n$'s for $n \geq 1$ are called of type $1$. Finally, a last step is added to only allow as isometries isomorphisms whose restriction to $Gx_0$ is a permutation which corresponds to the action of some $g \in G$ on $Gx_0$. This is technically more complicated and is obtained by adding cones at some points of $\overline{span} Gx_0$ which code the structure of $G$ and are called of type $2$.

The reader may get a geometric feeling of this proof by looking at the group $G=\{\pm Id,\pm R\}$ of $\R$-linear isometries on $\C$
where $R$ is the rotation of angle $\pi/2$. By adding cones on the unit ball at $\pm 1$ and $\pm i$, one allows the isometries in $G$ but also symmetries with respect to the axes. A way of correcting this is to add one well-placed smaller cone next to each element of $\{\pm 1, \pm i\}$ so that the only isometries in the new norm are those of $G$.

\begin{defi} Let $X$ be a real Banach space with norm $\|.\|$, let $G$ be a  group of isometries on $X$ such that $-Id \in G$, and let $(x_k)_{k \in K}$ be a possibly finite sequence of unit vectors of $X$. Let
$\Lambda=(\lambda_k)_{k \in K}$ be such that $1/2<\lambda_k<1$ for all
$k \in K$. The {\em $\Lambda,G$-pimple at $(x_k)_k$ for $\|.\|$} is the equivalent norm on $X$ defined by
$$\norm{y}_{\Lambda,G}=\inf\{\sum [[y_i]]_{\Lambda,G}: y=\sum y_i\},$$
where $[[y]]_{\Lambda,G}=\lambda_k \norm{y},$
whenever $y \in Vect(g.x_k)$ for some $k \in K$ and $g \in G$, and
$[[y]]_{\Lambda,G}=\norm{y}$ otherwise.
\end{defi}

 In other words, the unit ball
for $\norm{.}_{\Lambda,G}$ may be seen as the convexification of the union of the unit ball for $\norm{.}$ with line segments between
$gx_k/\lambda_k$ and $-gx_k/\lambda_k$ for each $k \in K$ and $g \in G$.

Some observations are in order. First of all
$(\inf_{k \in K}\lambda_k)\norm{.} \leq \norm{.}_{\Lambda,G} \leq \norm{.}.$ Any $g \in G$ remains an isometry in the norm $\norm{.}_{\Lambda,G}$. In \cite{Bel} Bellenot had defined the notion of $\lambda$-pimple at $x_0 \in X$, which corresponds to $(\lambda),\{-Id,Id\}$-pimple in our terminology.

Recall that a norm $\norm{.}$ is LUR at some point $x$, $\norm{x}=1$, if
$\forall 0<\epsilon \leq 2$, there exists $\lambda(x,\epsilon)<1$ such that whenever $\norm{y}=1$ and $\norm{x-y} \geq \epsilon$, it follows that $\norm{\frac{x+y}{2}} \leq \lambda(x,\epsilon)$.
Equivalently $\lim_n x_n=x$ whenever $\lim_n \norm{x_n}=\norm{x}$
and $\lim_n \norm{x+x_n}=2\norm{x}$. It is LUR when it is LUR at all points of the unit sphere. A norm is strictly convex if whenever $\norm{x}=\norm{y}=1$, the segment $[x,y]$ intersects the unit sphere in $x$ and $y$ only.
We recall a crucial result from \cite{Bel}.

\begin{prop}\label{bellenot}(Bellenot \cite{Bel}) Let $(X,\norm{.})$ be a real Banach space and let $\norm{x_0}=1$ so that
\begin{itemize}
\item (1) \norm{.} is LUR at $x_0$, and 
\item (2) there exists $\epsilon>0$ so that if $\norm{y}=1$ and
$\norm{x_0-y} < \epsilon$, then $y$ is an extremal point (i.e. an extremal point of the ball of radius $\norm{y}$).
\end{itemize}
Then given $\delta>0$, $B>0$ and $0<m<1$, there exists
a real $0<\lambda_0<1$ of the form 
$\lambda_0=\max(m,\lambda_0(\epsilon,\delta,B,\lambda(x_0,\eta(\epsilon,\delta,B))))<1$, so that whenever $\lambda_0 \leq \lambda <1$ and
$\norm{.}_{\lambda}$ is the $\lambda$-pimple at $x_0$, then 
\begin{itemize}
\item (3) $m\norm{.} \leq \norm{.}_{\lambda} \leq \norm{.}$,
\item (4) if $1=\norm{y}>\norm{y}_{\lambda}$ then $\norm{x_0-y} <\delta$
or $\norm{x_0+y}<\delta$,
\item (5) $x_{\lambda}=\lambda^{-1}x_0$ is the only isolated extremal point of $\norm{.}_{\lambda}$ which satisfies $\norm{x/\norm{x}-x_0} <\epsilon$,
\item (6) if $w$ is a vector so that $x_\lambda$ and $x_\lambda+w$ are endpoints of a maximal line segment in the unit sphere of $\norm{.}_\lambda$, then $B \geq \norm{w} \geq \lambda^{-1}-1$.
\end{itemize}
\end{prop}

For more details we refer to \cite{Bel}. We generalize this result to $(\Lambda,G)$-pimples in a natural manner which for the $\Lambda$ part is inspired from $\cite{Bel}$. Write
$\Lambda \leq \Lambda'$ to mean $\lambda_k \leq \lambda_k^{\prime}$ for all $k \in K$, if $\Lambda=(\lambda_k)_k$ and
$\Lambda'=(\lambda_k^{\prime})_k$.

\begin{prop}\label{Gbellenot}
Let $(X,\norm{.})$ be a real Banach space, let
$G$ be a group of isometries on $X$ containing $-Id$ and let $(x_k)_{k \in K}$ be a possibly finite sequence of unit vectors of $X$. Assume
\begin{itemize}
\item (1)' $\norm{.}$ is strictly convex on $X$ and LUR in $x_k$ for each $k \in K$, and 
\item (2)' for all $k \in K$, $c_k:=\inf\{\norm{x_j-gx_k}: j \in K, g \in G, (j,g) \neq (k,Id)\}>0$.
\end{itemize}
Then given  $\delta>0$, $B=(b_k)_k>0$ and $0<m<1$, there exists
$\Delta=(\delta_k)_k$ with $\delta_0 \leq \delta$ and for all $k \geq 1$,
$\delta_k \leq \min(\delta_{k-1},c_k/4,1-\lambda(x_k,c_k))$,
and $0<\Lambda_0=(\lambda_{0k})_k<1$ with for all $k$, $\lambda_{0k}=\max(m,\lambda_{0}^{\prime}(\epsilon_k,\delta_k,b_k,\lambda(x_0,\eta(\epsilon_k,\delta_k,b_k))))<1$, so that whenever $\Lambda_0 \leq \Lambda <1$ and
$\norm{.}_{\Lambda,G}$ is the $\Lambda,G$-pimple at $(x_k)_k$, then 
\begin{itemize}
\item (3)' $m\norm{.} \leq \norm{.}_{\Lambda,G} \leq \norm{.}$,
\item (4)' if $1=\norm{y}>\norm{y}_{\Lambda,G}$ then  $\exists g \in G, k \in K: \norm{gx_k-y} <\delta_k$
\item (5)' $x_{k,\lambda}=\lambda_k^{-1}x_k$ is the only isolated extremal point of $\norm{.}_{\Lambda,G}$ which satisfies $\norm{x/\norm{x}-x_k} <\epsilon_k$,
\item (6)' if $w$ is a vector so that $x_{k,\lambda}$ and $x_{k,\lambda}+w$ are endpoints of a maximal line segment in the unit sphere of $\norm{.}_{\Lambda,G}$, then $b_k \geq \norm{w} \geq \lambda_k^{-1}-1$.
\end{itemize}
\end{prop}

\pf Proposition \ref{bellenot} corresponds to the case $G=\{-Id,Id\}$ and $K$ a singleton. We shall deduce the general case from Proposition \ref{bellenot} and from the fact that for well-chosen $\Lambda$, the closed unit ball of the $\Lambda,G$-pimple at $(x_k)_k$ is equal to $B_0$, the union over $k \in K$ and $g \in G$ of the closed unit balls $B_{k,g}$  of the $\lambda_k$-pimples $\norm{.}_{\lambda_k,g}$ at $gx_k$. Let $B$ denote the closed unit ball for $\norm{.}$.

Note that by (1)', Proposition \ref{bellenot} (1)(2) apply in any $x_k, k \in K$, for any $\epsilon>0$. Let $\epsilon_k=c_k/2$.
Let  $\lambda_{0k} \geq \max(m,\lambda_{0}^{\prime}(\epsilon_k,\delta_k,b_k,\lambda(x_k,\eta(\epsilon_k,\delta_k,b_k))))$ given by Proposition \ref{bellenot} in $x_k$ for $\epsilon=\epsilon_k$, with $1-\lambda_{0k}^{-1} \leq c_k/6$ for all
$k \in K$ and with $\lim_{k \rightarrow +\infty} \lambda_{0k}=1$ if $K$ is infinite. The limit condition on $\lambda_{0k}$ ensures that $B_0$ is closed. Assuming $x,y \in B_0$ and $\frac{x+y}{2} \notin B_0$ let
$(k,g)$ and $(l,h)$ be such that $x \in B_k^g$ and $y \in B_l^h$. By convexity of $B_k^g$ and $B_l^h$, either $k \neq l$ (e.g. $k<l$),
 or $k=l$
and $g \neq \pm h$, and $x \in B_k^g \setminus B$,
$y \in B_l^h \setminus B$, i.e. 
$\norm{x}_{\lambda_k,g}<\norm{x}$ and
$\norm{y}_{\lambda_l,h}<\norm{y}$.
Therefore by (4) applied to $x$ for the $\lambda_k$-pimple at $gx_k$, and up to replacing $g$ by $-g$ if necessary,
$\norm{gx_k-x} < \delta_k$. Likewise
$\norm{hx_l-y} <\delta_l$.
Then $$\norm{\frac{x+y}{2}-\frac{gx_k+hx_l}{2}}<\frac{\delta_k+\delta_l}{2} \leq \delta_k.$$
Since $\norm{gx_k-hx_l} \geq c_k$ by (2)', it follows by LUR of $\norm{.}$ in $gx_k$ that
$$\norm{\frac{gx_k+hx_l}{2}} \leq \lambda(gx_k,c_k)=\lambda(x_k,c_k),$$
and
$$\norm{\frac{x+y}{2}} \leq \delta_k+\lambda(x_k,c_k)\leq 1,$$
a contradiction. Therefore $B_0$ is closed convex and $B_0$ is equal to the closed unit ball of the $\Lambda,G$-pimple at $(x_k)_k$.
Equivalently
$$\norm{.}_{\Lambda,G}=\inf_{k \in K,g \in G}\norm{.}_{\lambda_k,g}.$$
In fact, since whenever $x \in B_k^g \setminus B$ and
$y \in B_l^h \setminus B$ with $B_k^g \neq B_l^h$ and $k \leq l$, and up to replacing $g$ by $-g$ or $h$ by $-h$ if necessary, we have
$$\norm{x-y} \geq \norm{gx_k-hx_l}-\norm{x-gx_k}-\norm{y-hx_l} \geq
c_k-\delta_k-\delta_l \geq c_k/3,$$
it follows that for any $x$ such that $\norm{x}_{\Lambda,G}<\norm{x}$, there exists a unique $(g,\lambda_k)$ such that $\norm{x}_{\lambda_k,g}<\norm{x}$, and $\norm{x}_{\Lambda,G}=\norm{x}_{\lambda_k,g}$.

\

We now prove (3)'-(6)'. (3)' is obvious from (3) for each $\norm{.}_{\lambda_k,g}$.
For (4)' assume $1=\norm{y}>\norm{y}_{\Lambda,G}$ then as we have just observed, there exist $g,k$ such that
$1=\norm{y}>\norm{y}_{\lambda_k,g}$, so from (4),
$\norm{gx_k-y}<\delta_k$ or $\norm{-gx_k-y}<\delta_k$.

\

To prove (5)' note that if $\norm{x/\norm{x}-x_{k,\Lambda}}<\epsilon_k$ then whenever $g \neq Id$ or $k \neq l$,
$$\norm{x/\norm{x}-gx_l} > \norm{gx_l-x_k}-\norm{x_k-x_{k.\Lambda}}-\epsilon_k$$
$$\geq c_k-(1-\lambda_k^{-1}) -\epsilon_k \geq c_k/2-(1-\lambda_{k0}^{-1}) \geq \delta_k.$$
Therefore by (4)' $\norm{x}=\norm{x}_{\lambda_l,g}$ whenever
$g \neq Id$ or $k \neq l$,
and so 
$\norm{x}_{\Lambda,G}=\norm{x}_{\lambda_k}$.
Now if $x$ is an isolated extremal point of $\norm{.}_{\Lambda,G}$, it is therefore an isolated extremal point of $\norm{.}_{\lambda_k}$ and by (5), $x=x_{k,\Lambda}$.

\

The proof of (6)' is a little bit longer. Write
$S_k^g$ the unit sphere for $\norm{.}_{\lambda_k,g}$,
$S_{\Lambda}^G$ the unit sphere for $\norm{.}_{\Lambda,G}$,
$S$ the unit sphere for $\norm{.}$, $S'$ the set of points of $S$ on which $\norm{.}_{\Lambda,G}=\norm{.}$. As we know,
$S_{\Lambda}^G=S' \cup (\cup_{k,g} (S_k^g \setminus S))$.

As we noted before, whenever $x \in S_k^g \setminus S$,
$y \in S_l^h \setminus S$, with $S_k^g \neq S_l^h$, it follows that
$\norm{x-y} \geq c_{\min(k,l)}/3$. So for 
$x \in S_k^g \setminus S$, $\norm{x-y} \geq \frac{1}{3}\min\{c_i, i \leq k\}$ whenever $y$ belongs to some $S_l^h \setminus S$, with $S_k^g \neq S_l^h$.
 Therefore a line segment in $S_{\Lambda,G}$ containing points both in $S_k^g \setminus S$ and $S_l^h \setminus S$ with $S_k^g \neq S_l^h$ must have a subsegment included in $S$, but this contradicts the strict convexity of $\norm{.}$. 

We deduce that if $[x_{k,\Lambda},x_{k,\Lambda}+w]$ is a maximal line segment in $S_{\Lambda,G}$, it is a line segment in $S_k^{Id}$. It is now enough to prove that it cannot be extended in $S_k^{Id}$, then by (6) applied for
$\norm{.}_{\lambda_k,Id}$, $b_k \geq \norm{w} \geq \lambda_k^{-1}-1$.

But for any strict extension $[x_{k,\lambda},y]$ of
$[x_{k,\Lambda},x_{k,\Lambda}+w]$ in $S_k^{Id}$, either
$[x_{k,\lambda},y] \subset S_k^{Id} \setminus S \subset S_{\Lambda}^G$ and the maximality in $S_{\Lambda}^G$ is contradicted,
or there exists a sequence $(y_n)_n$ of distinct points converging
to $x_{k,\Lambda}+w$ in $[x_{k,\lambda},y]$ with $y_n \in S$ for all $n$,
but this again contradicts the strict convexity of $\norm{.}$. \pff

\begin{theo}\label{groupe} Let $X$ be a separable real Banach space with an LUR-norm $\norm{.}$ and let $G$ be a countable group of isometries on $X$ such that $-Id \in G$.
Assume that there exists a unit vector $x_0$ in $X$ which separates $G$ and such that the orbit $Gx_0$ is discrete. Then $X$ admits an equivalent norm $\||.\||$ such that $G$ is the group of isometries on $X$ for $\||.\||$. \end{theo}

\pf Since $Gx_0$ is discrete and $x_0$ separates $G$, let $\alpha \in ]0,1[$ be such that
$\norm{x_0-gx_0} \geq \alpha$, for all $g \neq Id$. 

Let $V_0=\overline{span}\{gx_0, g \in G\}$ and let $y_0=x_0$.
If $V_0 \neq X$ then it is possible to pick a possibly finite sequence 
$(y_n)_{n \geq 1}$ such that, if $V_n:=\overline{span}\{gy_k, k \leq n, g \in G\}$, we have that $y_n \notin V_{n-1}$ for all $n \geq 1$ and
$\cup_n V_n$ is dense in $X$.

\

Let $(u_n)_{n \geq 1}$ be a (possibly finite) enumeration of $\{gx_0, g \in G\setminus\{\pm Id\}\} \cup \{y_k, k \geq 1\}$. Then define a (possibly finite) sequence $(x_n)_n$ of unit vectors of $X$ by induction as follows.
Assume $x_0,\ldots,x_{n-1}$ are given.

If $u_n=y_k$ for some $k \geq 1$ then let $E=span(V_{k-1},y_k)$. Pick some $z_n \in E$ such that $\norm{z_n} \in [\alpha/10,\alpha/5]$ and 
$d(z_n,V_{k-1})=\alpha/10$, and let $x_n=a_n x_0+z_n$ where $a_n>0$ is such that $\norm{x_n}=1$. Such an $x_n$ will be called of type $1$.

If $u_n$ is of the form $gx_0$ then we shall pick some $\alpha_n \in [\alpha/10,\alpha/5]$ and define $z_n=\alpha_n gx_0$, $x_n=a_n x_0+z_n$ with $a_n>0$ and $\norm{x_n}=1$. Such an $x_n$ will be called of type $2$.
The choice of $\alpha_n$ will be made more precise later. Let us first observe a few facts.

By construction, $X$ is the closed linear span of $\{gx_n, g \in G, n\}$ (actually only $x_0$ and $x_n$'s of type $1$ are required for this).
Note that for all $n$, $\norm{z_n} \leq \alpha/5$ and therefore
$a_n \in [1-\alpha/5,1+\alpha/5]$; and obviously $x_0$ may also be written
$x_0=a_0x_0+z_0$ with these conditions.
We now evaluate $\norm{x_n-gx_m}$ for all $(n,Id) \neq (m,g)$.

If $g \neq Id$ then 
$\norm{x_n-gx_m} =\norm{a_n x_0+z_n-ga_m x_0-z_m}$ therefore
$$\norm{x_n-gx_m} \geq \norm{x_0-gx_0}-|1-a_n|-|1-a_m|-\norm{z_n}-\norm{z_m} \geq \alpha/5.$$
If $g=Id$, without loss of generality assume $n>m$. If $x_n$ is of type $1$ then, if $k$ is such that $x_n$ is associated to $y_k$, the vector $gx_m$ is in $V_{k-1}$ and $\norm{x_n-gx_m} \geq d(x_n,V_{k-1})=\alpha/10$. If $x_n$ is of type $2$ and $x_m$ of type $1$ then
$\norm{x_n-gx_m}=\norm{x_m-g^{-1}x_n} \geq d(x_m,V_0) \geq \alpha/10$.

\

It now remains to study the more delicate case when $x_n$ and $x_m$ both are of type $2$, or one is of type $2$ and the other is $x_0$. We describe how to choose the $x_n$'s of type $2$, i.e. how to choose each corresponding $\alpha_n \in [\alpha/10,\alpha/5]$
in the definition of $x_n$ to obtain good estimates for $\norm{x_n-x_m}$ in that case. This will be done by induction.
To simplify the notation, we shall denote $(x^{\prime}_n)_{n \in N}$ the subsequence $(x_{k_n})_{n \in N}$ corresponding to the $x_k$'s of type $2$, with $N=\{1,\ldots,|G|-2\}$ or $N=\N$ according to the cardinality of $G$, 
and we shall write $x^{\prime}_n=b_n x_0+\beta_n g_n x_0$, where $g_n$ is the associated element of $G\setminus \{\pm Id\}$, $b_n=a_{k_n}$ and
$\beta_n=\alpha_{k_n}$. Write $x^{\prime}_0=x_0$.

 Let $\forall m \geq 1, I_m^0=[\alpha/10,\alpha/5]$.
For $\beta \in [\alpha/10,\alpha/5]$, let $x^{\prime}_m(\beta)=b_m(\beta)x_0 +\beta g_mx_0$ where $b_m(\beta)>0$ is such that $\norm{x^{\prime}_m(\beta)}=1$. 

We observe that $\norm{x^{\prime}_m(\beta)-x^{\prime}_m(\gamma)} \geq \frac{\alpha}{2}|\beta-\gamma|$.
Indeed if $x^{\prime}_m(\beta)-x^{\prime}_m(\gamma)=(\beta-\gamma)\epsilon$ with $\norm{\epsilon} <\alpha/2$ and $\beta \neq \gamma$, then
$$(b_m(\beta)-b_m(\gamma))x_0=(\gamma-\beta)(g_m x_0-\epsilon),$$
so $g_mx_0-\epsilon=\pm \norm{g_m x_0-\epsilon} x_0$. 
If for example $\pm=-$ in this equality, then
$$\norm{g_m x_0+x_0}=\norm{\epsilon+(1-\norm{g_m x_0-\epsilon})x_0} \leq 2\norm{\epsilon} <\alpha,$$
and by separation, $g_m=-Id$, a contradiction. Similarly the case $\pm=+$ would imply $g_m=Id$.

Now for all $m \geq 1$ divide $I_m^0=[\alpha/10,\alpha/5]$ in three successive intervals of equal length $\alpha/30$. Since
$$\norm{x^{\prime}_m(\beta)-x^{\prime}_m(\gamma)} \geq \frac{\alpha}{2}|\beta-\gamma|\geq \frac{\alpha^2}{60}$$
whenever $\beta$ is in the first and $\gamma$ in the last interval, it follows that there exists an interval $I_m^1 \subset I_m^0$ of length
$\alpha/30$ (which is either the first or the last subinterval), such that
$$\beta \in I_m^1 \Rightarrow \norm{x^{\prime}_m(\beta)-x^{\prime}_0} \geq \frac{\alpha^2}{120}.$$
We then pick  $\beta_1$ in $I_1^1$ and fix $x^{\prime}_1=x^{\prime}_1(\beta_1)$. Therefore we have ensured
$$\norm{x^{\prime}_1-x^{\prime}_0} \geq \frac{\alpha^2}{120}.$$
Assume selected $\beta_1,\ldots,\beta_{n-1}$, $x^{\prime}_1,\ldots,x^{\prime}_{n-1}$ associated, and for $0 \leq i \leq n-1$ and
 $m \geq i$, decreasing in $i$ intervals $I_m^i$ of length $\frac{\alpha}{10.3^i}$.
For any $m \geq n-1$, dividing $I_m^{n-1}$ in three subintervals and picking the first or the last, we find by the same reasoning as above
$I_m^n \subset I_m^{n-1}$ of length  $\frac{\alpha}{10.3^n}$ with
 $$\beta \in I_m^n \Rightarrow \norm{x^{\prime}_m(\beta)-x^{\prime}_{n-1}} \geq \frac{\alpha^2}{40.3^n}.$$
We then pick  $\beta_n$ in $I_n^n$ and fix $x^{\prime}_n=x^{\prime}_n(\beta_n)$. Therefore for all $k<n$, $\beta_n \in I_n^n \subset I_n^{k+1}$ and we have ensured
$$\forall 0 \leq k<n, \norm{x^{\prime}_n-x^{\prime}_k} \geq \frac{\alpha^2}{40.3^{k+1}}.$$

We have
finally proved that for all $k$, $$\inf\{\norm{x_n-gx_k}, n \geq k, g \in G, (n,g) \neq (k,Id)\} \geq \frac{\alpha^2}{40.3^{k+1}},$$
and so
$$\inf\{\norm{x_n-gx_k}, n, g \in G, (n,g) \neq (k,Id)\} \geq \frac{\alpha^2}{40.3^{k+1}},$$
 therefore (2)' in Proposition \ref{Gbellenot} is satisfied; and (1)' is clearly satisfied since $\norm{.}$ is LUR.

We then define $\||.\||$ as the $\Lambda,G$-pimple at $(x_n)_n$ for
$\Lambda=(\lambda_n)_n$ associated 
to $\epsilon_n,b_n$ so that Proposition \ref{Gbellenot} applies and such that $b_n>b_{n+1}$, $1/2 \leq \lambda_n<\lambda_{n+1}$ and $\lambda_n^{-1}-1>2b_{n+1}$ for all $n$. This is possible by induction and the expression of $\Lambda_0$ in Proposition \ref{Gbellenot}. 

\

Observe that $E=\{gx_n/\lambda_n, g \in G, n\}$ is the set of isolated extremal points of $\||.\||$. Indeed
for a point $x$ of $S_{\Lambda,G}$ either $\norm{x/\norm{x}-gx_k}<\epsilon_k$ for some $g,k$, in which case by (5)' $x=\lambda_k^{-1}gx_k$ if it is an isolated extremal point; or
$\norm{x/\norm{x}-gx_k} \geq \epsilon_k >\delta_k$ for all $g,k$ then by
(4)' $\norm{.}=\||.\||$ in a neighborhood of $x$ and then $x$ is not an isolated extremal point since $\norm{.}$ is LUR at $x$.

Therefore any isometry $T$ for $\||.\||$ maps $E$ onto itself. If $n<m$, $g \in G$, then $T$ cannot map $\lambda_n^{-1}x_n$ to
$\lambda_m^{-1} gx_m$. Indeed if $w$ (resp. $w'$) is a vector so that
$\lambda_n^{-1}x_n$ and $\lambda_n^{-1}x_n+w$ (resp.
$\lambda_m^{-1}gx_m$ and $\lambda_m^{-1}gx_m+w'$) are endpoints of a maximal line segment in the unit sphere of $\||.\||$, then since $g$ is an isometry for $\||.\||$ we may assume $g=Id$, and then by (6)',
$$\||w\|| \geq \frac{1}{2} \norm{w} \geq\frac{1}{2}(\lambda_n^{-1}-1)>b_{n+1} \geq b_m \geq \norm{w'} \geq \||w'\||.$$
It follows that for each $n$, the orbit $Gx_n$ is preserved by $T$.

\

We finally prove that $T$ belongs necessarily to $G$. Without loss of generality we may assume that $Tx_0=x_0$ and then by density it is enough to prove that $Tgx_n=gx_n$ for all $g \in G$ and any $x_n$ of type $1$ or equal to $x_0$.

Let $g \in G$, $g \neq \pm Id$. Let $x'$ be the associated vector of type $2$ of the form $x'=ax_0+\beta gx_0$. Then $$Tx'=ax_0+\beta Tgx_0=h(ax_0+\beta gx_0)$$ for some $h \in G$.
So $|a|\norm{x_0-hx_0}=\beta\norm{Tgx_0-hgx_0}$ and
$$\norm{x_0-hx_0} \leq \frac{\alpha/5}{1-\alpha/5}(1+2\||Tgx_0\||)
\leq \frac{\alpha}{4}(1+2\||x_0\||) <\alpha,$$
therefore by separation $h=Id$. It follows immediately that
$$Tgx_0=gx_0.$$
and this holds for any $g \in G$. Finally if $x_n$ is of type $1$, and $g \in G$, then
$$Tgx_n=T(a_ngx_0+gz_n)=a_ngx_0+Tgz_n,$$
and since $T(gx_n)$ is of the form $hx_n$ for some $h \in G$,
$$Tgx_n=a_nhx_0+hz_n.$$
Therefore $a_n\norm{gx_0-hx_0}=\norm{Tgz_n-hz_n}$ and by similar computations as above,
$$\norm{gx_0-hx_0} \leq \frac{3\alpha/5}{1-\alpha/5}<\alpha,$$ whence again by separation $g=h$ and 
$$Tgx_n=gx_n.$$ \pff

\subsection{Representable groups of linear isomorphisms}

In this subsection, we give sufficient conditions for a group of isomorphims on a Banach space $X$ to be representable in $X$.

\begin{theo}\label{avecRNP} Let $X$ be a separable real Banach space with the Radon-Nikodym Property and $G$ be a countable bounded group of isomorphisms on $X$, containing $-Id$, and such that some point separates $G$ and has discrete orbit. Then $X$ admits an equivalent norm for which $G$ is the group of isometries on $X$.
\end{theo}

\pf We may assume that every $g$ in $G$ is an isometry on $X$ by using the equivalent norm $\sup_{g \in G}\norm{gx}$. Then by a result of G. Lancien, \cite{L} Theorem 2.1, and since $X$ is separable and has the RNP, it may be renormed with an LUR norm without diminishing the group of isometries (this last fact is not written explicitely in \cite{L} but is obvious from the definition of the renorming). We are then in position to apply Theorem \ref{groupe}. \pff 

\begin{theo}\label{renormingG} Let $X$ be a separable real Banach space and $G$ be a finite
group of isomorphisms such that $-Id \in G$. Then $X$ admits an equivalent norm for which $G$ is the group of isometries on $X$.
\end{theo}

\pf By a classical theorem of Kadec (we refer to \cite{DGZ} about LUR-renorming questions) we may assume that the norm $\norm{.}$ on $X$ is LUR. Then we define an equivalent norm $\norm{.}_G$ on $X$ by
$$\norm{x}_G=(\sum_{g \in G}\norm{gx}^2)^{1/2}.$$
Since this is the $l_2$-sum of the LUR norm $\norm{.}$ with an equivalent norm, it is classical to check that it is also LUR, see \cite{DGZ} Fact 2.3, and obviously any $g \in G$ becomes an isometry for $\norm{.}_G$.
To apply Theorem
\ref{groupe} it therefore only remains to
find some $x_0$ such that $x_0 \neq gx_0$ for all $g \neq Id$. But if such an $x_0$ didn't exist then $Ker(Id-g)$ would have non-empty interior for some $g \neq Id$, but by linearity this would actually imply that $g=Id$.
\pff

\

Note that the condition in Theorem \ref{avecRNP} that some point separates $G$ and has discrete orbit implies directly that $G$ is closed (and discrete) in the strong operator topology and therefore also for the convergence of $T$ and $T^{-1}$ in the strong operator topology.
Conversely to Theorem \ref{avecRNP}:

\begin{lemm}\label{discrete} Let $X$ be a separable real Banach space and $G$ be a group of isomorphisms which is the group of isometries in some equivalent norm on $X$. If $G$ is countable then $G$ is discrete for the convergence of $T$ and $T^{-1}$ in the strong operator topology, and $G$ admits a separating point. If $X$ is finite dimensional and $G$ is countable then $G$ is finite. \end{lemm}

\pf The existence of a separating point is a consequence of the Theorem of Baire. Indeed for any $g \in G$, $g \neq Id$, the set of points which separate $g$ from $Id$, i.e. the set $X \setminus Ker(g-Id)$, is dense open, therefore the set of separating points is a $G_{\delta}$ dense set.

To prove that $G$ is discrete we may assume that the norm is such that $G$ is the group of isometries on $X$. It is classical to check that $G$ is Polish.
Indeed
since $X$ is separable, the unit ball $L_1(X)$ of $L(X)$ with the (relative) strong topology is Polish \cite{Kec}, page 14.  We define $\phi: G \to L_{1}(X) \times L_{1}(X)$ by $\phi(T)=(T, T^{-1})$ and note that
$\phi(G)$ is closed in 
$L_{1}(X) \times L_{1}(X)$ (this follows immediately from the fact that if $(T_{n})_{n \in \N}$ converges to $T$ in $L_{1}(X)$ and $(L_{n})_{n \in \N}$ converges to $L$ in $L_{1}(X)$, then $T_{n} L_{n}$ converges to $TL$ in $L_{1}(X)$).
Hence  $\phi(G)$ is a Polish space, and as $\phi$ is a bijection onto the image, $G$ is a Polish space with the induced topology by $\phi$.
We then conclude using the fact that every countable Polish group is a discrete space.
 Indeed if $G$ is a countable Polish group, then by \cite{Kec}, Theorem 6.2, $G$ is not a perfect space, that is, $G$ has an isolated point, therefore by the group property all points are isolated. 

Finally if $X$ is finite dimensional then the strong topology on $L_1(X)$ coincides with the usual one for which $L_1(X)$ is compact. So $\phi(G)$ as a discrete subset of $L_1(X) \times L_1(X)$ is therefore finite.\pff

\

Note however that it seems to remain unknown whether a group $G$ of isomorphisms, which is the group of isometries on a real Banach space $X$ in some equivalent norm, and which is countable, must have some separating point with discrete orbit. 

\

It is also natural to wonder whether the only role of the separation and discrete orbit hypothesis in Theorem \ref{avecRNP} is to guarantee the closedness of the group $G$. That is, for groups which are not closed for the convergence of $T$ and $T^{-1}$ in the strong operator topology, one may wish to generalize Theorem \ref{avecRNP} by showing that whenever $X$ is separable with the RNP, and
$G$ is a countable bounded group of isomorphisms containing $-Id$, then there exists an equivalent norm for which the group of isometries is equal to the corresponding closure $\overline{G}^{op}$. This however is false as proved by the next example.

\begin{exam}\label{exampledeC} Let $G$ be the group of rational rotations on $\C$. Then $\C$ cannot be renormed, as a real space, so that the group of $\R$-linear isometries on $X$ is $\overline{G}^{op}$. \end{exam}

\pf The set $\overline{G}^{op}$ is the set of rotations on $\C$. If $\||.\||$ is a new real norm on $\C$ which is invariant by rotations, then it is a multiple of the modulus. But then the symmetry with respect to the real axis is an isometry on $\C$ which does not belong to $\overline{G}^{op}$. \pff

\

However the next question remains open in general (i.e. for a space which does not have the RNP):

\begin{ques} Let $X$ be a separable real Banach space and let $G$ be an infinite countable bounded group of isomorphisms on $X$ such that $-Id \in G$, and some point separates $G$ and has discrete orbit. Does $X$ admit an equivalent norm for which $G$ is the group of isometries on $X$?
\end{ques}

\subsection{Representation of countable groups in Banach spaces}

 Jarosz conjectured that any group of the form $\{-1,1\} \times G$ (or $C \times G$ in the complex case) could be represented in any Banach space $X$ provided $\dim X \geq |G|$. From Theorem \ref{renormingG} and Theorem \ref{groupe} we obtain rather general answers to his question for countable groups and separable real spaces.

\begin{theo}\label{fini} Let $G$ be a finite group and $X$ be a separable real Banach space such that $\dim X \geq |G|$. Then $\{-1,1\} \times G$ is representable in $X$.
\end{theo}

\pf The group $\{-1,1\} \times G$ may be canonically represented as a group of isometries on $\ell_2(G)$: denoting
$(e_g)_{g \in G}$ the canonical basis of $\ell_2(G)$, associate to any
$(\epsilon,g)$ in $\{-1,1\} \times G$ the isometry $T_{\epsilon,g}$ defined on $\ell_2(G)$ by
$$T_{\epsilon,g}(\sum_{h \in G} \lambda_h e_h)=\epsilon \sum_{h \in G} \lambda_h e_{gh}.$$
Since $\dim X \geq |G|$, the space $X$ is isomorphic to the $l_2$ direct sum $l_2(G) \oplus_2 Y$, for some space $Y$. By associating to any
$(\epsilon,g)$ in $\{-1,1\} \times G$ the isometry $A_{\epsilon,g}$ defined on 
$\ell_2(G) \oplus_2 Y$ by
$$A_{\epsilon,g}(t,y)=(T_{\epsilon,g}(t),\epsilon y),$$
we see that $\{-1,1\} \times G$ is isomorphic to a group of isometries on $\ell_2(G) \oplus_2 Y$ containing $-Id$. Therefore Theorem
\ref{renormingG} applies 
to deduce that $\{-1,1\} \times G$ is isomorphic to the group of isometries on $X$ in some equivalent norm. \pff

\

By Lemma \ref{discrete} an infinite countable group is representable in a real space $X$ only if $X$ is infinite dimensional.
For finite groups, it seems to remain open whether the condition on the dimension is necessary in Theorem \ref{fini}. This is not the case when $|G|$ is an odd prime. Indeed, letting $p=|G|$, $G$ is then isomorphic to $\Z/p\Z$ and so $\{-1,1\}\times G$ is isomorphic to $\Z/2p\Z$ and therefore may be represented as the group $\{e^{i k \pi/p}Id, 0 \leq k \leq 2p-1\}$ of isometries on $\C$; so $\{-1,1\} \times G$ is representable in $\R^2$, and by the proof of Theorem \ref{fini}, in any separable real space of dimension at least $2$.
For other values of $|G|$ the question remains open:

\begin{ques} For arbitrarily large $n \in \N$, does there exist a group $G$ with $|G|=n$, such that $\{-1,1\} \times G$ is representable in a separable real Banach space $X$ if and only if $\dim X \geq n$?
\end{ques}

\

A group which is representable in a Banach space necessarily admits a normal subgroup with two elements.
Recall that reciprocally  any (resp. countable) group which admits a normal subgroup with two elements is representable in a (resp. the separable) Hilbert space \cite{St}. The next theorem shows that this extends to a wide class of spaces, including the classical spaces $c_0$, $C([0,1])$, $\ell_p, 1 \leq p <+\infty$, and $L_p, 1 \leq p<+\infty$.

A Banach space is said to have a symmetric decomposition if it is of the form $(\sum \oplus Y)_S$, for some space $S$ with a symmetric basis $(s_n)_n$, i.e. the norm on $X$ is given by
$\norm{(y_n)_n}=\norm{\sum_n \norm{y_n}s_n}$.
 
\

\begin{theo} \label{classical}
Let $G$ be a countable group which admits a normal subgroup with two elements and $X$ be an infinite-dimensional separable real Banach space with a symmetric decomposition which either is isomorphic to $c_0(Y)$ or to $l_p(Y)$ for some $Y$ and $1 \leq p <+\infty$, or has the Radon-Nikodym Property. Then $G$ is representable in $X$.
\end{theo}

\pf We first assume that $G$ is infinite. Let $\{1,j\}$ be a normal subgroup of $G$ with two elements, therefore $j$ commutes with any element of $G$. Let $G'$ be a subset of $G$ containing $1$ and such that $G=G' \cup jG'$ and $G'\cap jG'=\emptyset$. For $g \in G$ let $\epsilon_g=1$ if $g \in G'$ and $\epsilon_g=-1$ otherwise, and let
$|g|$ denote the unique element of $\{g,jg\} \cap G'$.

Write $X=(\sum \oplus Y)_S$ and index the decomposition on $G'$, i.e write an element of $X$ as $(y_g)_{g \in G'}$.
We associate to any
$g$ in $G$ the isomorphism $T_{g}$ defined on 
$X$ by
$$T_{g}((y_h)_{h \in G'})=(\epsilon_{g^{-1}h}y_{|g^{-1}h|})_{h \in G'}.$$
Observe that if $g,k \in G$, then
$$T_k T_g ((y_h)_h)=T_k((\epsilon_{g^{-1}h}y_{|g^{-1}h|})_h)
=(\epsilon_{k^{-1}h} \epsilon_{g^{-1}|k^{-1}h|} y_{|g^{-1}|k^{-1}h||})_h.$$
Since $j$ commutes with any element of $G$, we have
$|g^{-1}|k^{-1}h||=|g^{-1}k^{-1}h|$ and it is easy to see that
$\epsilon_{k^{-1}h} \epsilon_{g^{-1}|k^{-1}h|}=\epsilon_{g^{-1}k^{-1}h}$,
therefore
$$T_k T_g ((y_h)_h)= (\epsilon_{(kg)^{-1}h} y_{(kg)^{-1}h})_h=T_{kg}((y_h)_h).$$
 From this we deduce that the map $g \mapsto T_g$ is a group homomorphism, and therefore
we may assume that $G$ is a bounded group of isomorphisms on $X$ containing $-Id$ (here identified with $j$).

Let $x_0$ be a unit vector in the summand of the decomposition indexed by $1$.
We observe that $\norm{x_0-(-x_0)}=2$ and that for any $g \in G, g \notin\{-Id,Id\}$,
$$\norm{x_0-gx_0} \geq c,$$
where $c$ is the constant of the basis $(s_g)_{g \in G'}$ of $S$.
Therefore $x_0$ separates $G$ and has discrete orbit. Finally, when $X$ has the RNP, Theorem \ref{avecRNP} applies.

When $X$ is isomorphic to $c_0(Y)$ or $\ell_p(Y)$ for some $1 \leq p <+\infty$, we use the existence of a LUR norm on $X$ for which the $T_g$'s are isometries. The existence of the LUR norm may be found in the Appendix, Lemma \ref{appendix} for $\ell_p(Y)$, Lemma \ref{appendixbis} for $c_0(Y)$, modulo the result of Kadec that any separable space $Y$ has an equivalent LUR norm. Therefore $G$ is representable as a group of isometries containing $-Id$ for an LUR norm on $X$. Any unit vector $x_0$ in the first summand of the decomposition separates $G$ and has discrete orbit, therefore  Theorem \ref{groupe} applies.

Finally in the case when $G$ is finite, we may index a symmetric decomposition of
$X$ on $\cup_{i \in \N} G_i^{\prime}$ where the $G_i^{\prime}$ are disjoint copies of $G'$. We may then use the previous method to represent $G$, up to renorming, as a group of isometries containing $-Id$ on each space spanned by the sum of the summands indexed on $G_i^{\prime}$, and therefore globally as
a group of isometries containing $-Id$ on X. The rest of the proof is as before. \pff

\begin{coro}\label{Lp} A countable group is representable in the real space $c_0$, resp. $C([0,1])$, $l_p$ for $1 \leq p<+\infty$, $L_p$ for $1 \leq p<+\infty$, if and only if it admits a normal subgroup with two elements. \end{coro}

From Theorem \ref{classical} we may also deduce the following theorem.

\begin{theo} \label{deno}
Let $G$ be a countable group and $X$ be an infinite-dimensional separable real Banach space which contains a complemented subspace with a symmetric basis. Then $\{-1,1\} \times G$ is representable in $X$.
\end{theo}

\pf By Theorem \ref{fini} we may assume that $G$ is infinite. Let $Y$ be a complemented subspace $Y$ of $X$ with a symmetric basis, and 
write $X=Y \oplus Z$. Since a symmetric basis is unconditional, $Y$ is either reflexive or contains a complemented subspace isomorphic to $c_0$ or $l_1$, therefore we may assume that $Y$ has the Radon-Nikodym Property or is isomorphic to $c_0$. By Theorem \ref{classical}
we may assume that $\{-1,1\} \times G$ is a group of isometries on $Y$ containing $-Id$ (here identified with $(-1,1_G)$).

When $Y$ has the RNP, we may by applying the result of Lancien \cite{L} Theorem 2.1 also assume that the new norm is LUR. Since $Z$ is separable we may also assume it is equipped with an LUR norm, and we equip $X$ with the $l_2$-sum norm $\||.\||$, i.e. $X=Y \oplus_2 Z$. It is classical that the norm $\||.\||$ is LUR on $X$.

Furthermore, for any $(\epsilon,g)$ in $\{-1,1\} \times G$, the map $A_{\epsilon,g}$ defined on 
$X= Y \oplus_2 Z$ by
$$A_{\epsilon,g}(y,z)=((\epsilon,g).y,\epsilon z)$$
is an isometry on $X$ for $\||.\||$. Therefore $\{-1,1\} \times G$ is isomorphic to a group of isometries on $(X,\||.\||)$ containing $-Id$.
As in the proof of Theorem \ref{classical}, the point $x_0=e_1$
separates $G$ and has discrete orbit, where $e_1$ is the first vector of the symmetric basis of $Y$, so finally Theorem
\ref{groupe} applies.

When $Y$ is isomorphic to $c_0$, we may use Lemma \ref{appendixbis} to see $\{-1,1\} \times G$ as a group of isometries containing $-Id$ for an LUR norm on $Y$. The rest of the proof is as in the first case. \pff

\

Observe that Theorem \ref{deno} applies whenever $X$ is a subspace of $\ell_p$, $1 \leq p <+\infty$, or, by Sobczyk's Theorem, \cite{LT} 
Th. 2.f.5, whenever $X$ is separable and contains a copy of $c_0$.

\

Because of Theorem \ref{classical}, it is natural to ask whether Theorem \ref{fini} and Theorem \ref{deno} extend to the case when one replaces groups of the form $\{-1,1\} \times G$
by  groups which admit a normal subgroup with two elements.
We provide  examples to show that the answer is negative in general. 

The space denoted $X_{GM}$ is the real HI space of W.T. Gowers and B. Maurey
\cite{GM1}. Every operator on $X_{GM}$ is of the form $\lambda Id+S$, $\lambda \in \R$, $S$ strictly singular, and therefore every isometry is of the form $\pm Id+S$ (actually we shall see in the last section that every isometry on $X_{GM}$ is of the form $\pm Id+K$, $K$ compact, but this refinement will not be needed here). The complex version of $X_{GM}$ is such that every isometry is of the form $\lambda Id+S$, $\lambda \in C$, $S$ strictly singular (actually, $S$ compact, by \cite{RR}).

\begin{prop}\label{HIreel} Any group which is representable in the real (resp. the complex) $X_{GM}$ is of the form $\{-1,1\} \times G$ (resp. $C \times G$). In particular a finite group is representable in the real $X_{GM}$ if and only if it is of the form $\{-1,1\} \times G$. \end{prop}

\pf The last part of the proposition is a consequence of the initial part and of Theorem \ref{fini}. We prove the initial part. Let $H$ be the group of isometries on the real (resp. complex) $X_{GM}$ in some equivalent norm. 
Let $G$ be the subgroup of $H$ of isometries of the form $Id+S$, $S$ strictly singular.
For $T \in H$, let $\lambda_T$ be the element of $\{-1,1\}$ (resp. $C$) such that $T-\lambda_T Id$ is strictly singular. It is then easy to see, using the ideal properties of strictly singular operators, that by 
mapping $T$ to $(\lambda_T,T/\lambda_T)$ we provide an isomorphism of $H$ onto the group $\{-1,1\} \times G$ (resp. $C \times G$). \pff

\begin{prop}\label{symmetricplusHI} Let $S$ be a Banach space with a symmetric basis. Any group which is representable in $S \oplus X_{GM}$ is of the form $\{-1,1\} \times G$ in the real case (resp. $C \times G$ in the complex case). In particular, in the real case, a countable group is representable in  $S \oplus X_{GM}$ if and only if it is of the form $\{-1,1\} \times G$. \end{prop}

\pf The last part of the proposition is a consequence of the initial part and of Theorem \ref{deno}. We prove the initial part. 
Let $X=S \oplus X_{GM}$. We observe that, since $S$ and $X_{GM}$ are totally incomparable, any operator $T$ on $X$ may be written as a matrix 
of the form $$\begin{pmatrix} A & s_1 \\ s_2 & \lambda_T Id+s \end{pmatrix},$$
where $A \in L(S)$, and $s_1 \in L(X_{GM},S), s_2 \in L(S,X_{GM}), s \in L(X_{GM})$ are strictly singular; and $\lambda_T \neq 0$ if $T$ is an isomorphism. If $T$ is an isometry then since $T_{|X_{GM}}$ is a strictly singular perturbation of $\lambda_T i_{X_{GM},X}$,
where $i_{X_{GM},X}$ denotes the canonical
 injection of $X_{GM}$ into $X$, $\lambda_T$ must belong to $\{-1,1\}$ (resp. $C$).

Let $H$ be the group of isometries on $S \oplus X_{GM}$ for some equivalent norm. 
Let $G$ be the subgroup of $H$ defined by $G=\{T \in H: \lambda_T=1\}$. Clearly
mapping $T$ to $(\lambda_T,T/\lambda_T)$ we provide an isomorphism of $H$ onto the group $\{-1,1\} \times G$ (resp. $C \times G$). \pff

\

It remains open for a given separable infinite dimensional real space $X$ exactly which finite (resp. countable) groups are representable. We have the maximal case of $c_0$, $C([0,1])$, $\ell_p, 1 \leq p <+\infty$ or $L_p, 1 \leq p<+\infty$, in which all countable groups admitting a normal subgroup with two elements are representable, and the minimal case of $X_{GM}$, in which only groups of the form $\{-1,1\} \times G$ are representable. Apparently quite various situations may occur. Indeed we also show that a space constructed in \cite{F}  provide a third example which is "in between" the cases of $\ell_p$ and $X_{GM}$:
in the following $X(\C)$ denotes, seen as real, the separable complex space defined in 
\cite{F} on which every $\R$-linear operator is of the form $\lambda Id+S$, where $\lambda \in \C$ and $S$ is strictly singular.

\begin{prop}\label{HIcomplexe} The class of finite groups representable in $X(\C)$ is neither equal to the class of finite groups which admit a normal subgroup with two elements, nor to the class of finite groups of the form $\{-1,1\} \times G$.
\end{prop}

\pf 
For any $n \in \N, n \geq 1$, the group $\{e^{i k \pi/2n}Id, 0 \leq k \leq 4n-1\} \simeq \Z/4n\Z$ is a finite group of isomorphisms on $X(\C)$ containing $-Id$. Therefore by Theo\-rem \ref{renormingG} it is representable in $X(\C)$; however it is not of the form $\{-1,1\} \times G$.

On the other hand,
let $\{1,i,j,k\}$ be the generators of the algebra $\HH$ of quaternions, and let $G$ be the group $\{\pm 1,\pm i, \pm j, \pm k\}$. The
 group $\{-1,1\}$ is a normal subgroup of $G$ with two elements, and we prove that $G$ is not representable in $X(\C)$.

Assume on the contrary that $\alpha$ is an isomorphism from $G$ onto $H$,
where
$H$ is the group of isometries on $X(\C)$ in some equivalent norm.
 Since $-Id \in H$, $(-Id)^2=Id$ and $-1$ is the only element of square $1$ in $G\setminus\{1\}$, we have $\alpha(-1)=-Id$. 
Therefore from $ij=-ji$ we deduce $\alpha(i)\alpha(j)=-\alpha(j)\alpha(i)$. Let, for $T$ an operator on $X(\C)$, $\lambda_T$ be the unique complex number such that $T-\lambda_T Id$ is strictly singular. The map $T \mapsto \lambda_T$ induces an homomorphism of $H$ into $C$.
 We deduce $\lambda_{\alpha(i)}\lambda_{\alpha(j)}=-\lambda_{\alpha(j)}\lambda_{\alpha(i)}$, which is impossible in $C$. \pff

 \section{Complex structures up to isometry}

In this section we study complex structures on real Banach spaces up to isometry. This is related to the problem of representation of groups in Banach spaces, as we shall study the representation of the circle group $C$ and of the group of $\R$-linear isometries of $\C$ in real Banach spaces, to obtain uniqueness or non-uniqueness properties of isometric complex structures on a Banach space according to choices of equivalent renorming.

Note that it is immediate that two complex structures (in the isometric sense) $X^I$ and $X^J$ on a real Banach
space $X$ are isometric if and only if $I$ and $J$ are isometrically
conjugate, i.e. there exists an $\R$-linear isometry $P$ on $X$ such that $J=PIP^{-1}$ (the
operator $P$ is then a $\C$-linear isometry from $X^I$ onto $X^J$).

\subsection{The classical case of $\ell_2$}

Recall that the space $\ell_2$ admits a canonical complex structure associated to the isometry $J$ on $\ell_2=\ell_2 \oplus_2 \ell_2$ defined by $J(x,y)=(-y,x)$, i.e. this complex structure is the complexification $\ell_2 \oplus_{\C}^2 \ell_2$ of $\ell_2$.

 \begin{prop}\label{l2} The space $\ell_2$ admits a unique structure up to isometry.
\end{prop} 

\pf It is enough to prove that whenever $A$ is an isometry
 on  $\ell_2$ satisfying $A^{2}=-Id$,
$\ell_2^A$ is $\C$-linearly isometric to the canonical complex structure on $\ell_2$, i.e.
 there exists an orthonormal basis
 $(u_n)_{n \in \N}$ of $\ell_2$ such that decomposing $\ell_{2}=[(u_{2n-1})_{n \in \N}] \oplus [(u_{2n})_{n \in \N}]$  the  matrix  of $A$ is
 
 $$A=\left(\begin{array}{ll}0 & -Id\\
 \\Id &  \ \ \ 0 \end{array}\right).$$

We note the following fact. Fix a non-zero vector $x$ in
$\ell_2$. Since $A^2 =-Id$,
 the subspace $[x,Ax]$ generated by $x$ and $Ax$  is invariant
 by $A$ and has dimension $2$. Take an orthonormal basis $\{u, v \}$ of
 $[x, Ax]$. Then  the restriction of $A$  to $[u, v]$ is a rotation of
 angle $\theta$ for some $\theta \in \R$.  Moreover, $\theta= \pi/2$ or
 $\theta = 3 \pi /2$, because $A^2=-Id$. Therefore $Au=v$ and $Av=-u$
 or $Au=-v$ and $Av=u$. In particular, there exists an orthonormal
 basis $\{u_{1}, u_{2} \}$ of $[x, Ax]$  such that $Au_1=u_2$
 and $Au_2=-u_1$.

We also observe that since the adjoint operator of $A$ is $-A$,
 the orthogonal subspace  $[u_1, u_2]^{\perp}$  of $[u_1,  u_2]$
 is also invariant by $A$.

Let now $(e_n)_{n \in \N}$ be an orthonormal  basis of $\ell_2$. 
 By the fact there
 exists an orthonormal basis $\{u_1, u_2 \}$ of  $X_{1}=[e_1, Ae_1]$
  such that $Au_1=u_2$ and $Au_2=-u_1$. Pick $n_1$ the smallest
 $i \in \N$ verifying $e_{i} \not \in X_1$. Consequently $2 \leq n_1$.

We know that $\ell_{2}=X_{1} \oplus X_1^{\perp}$. So there exists some
 $f_2 \in X_1^{\perp}$ satisfying $e_{n_1} - f_2 \in X_1$. Denote
 $X_{2}=[f_2, Af_2]$. By our observation  $X_{2} \subset X_1^{\perp}$.
 Now by restricting $A$ to $X_2$, again by the fact there exists an orthogonal
 basis $\{u_3, u_4 \}$ of  $X_{2}$ such that $Au_3=u_4$ and $Au_4=-u_3$.
 Fix $n_2$ the smallest $i \in \N$ verifying $e_{i} \not \in X_1 \oplus X_2$.
 Thus $3 \leq n_2$.

Since $\ell_{2}=X_1 \oplus X_2 \oplus (X_1 \oplus X_2)^{\perp}$, there exists $f_{3} \in (X_1 \oplus X_2)^{\perp}$ verifying $e_{n_2}-f_3 \not \in X_1 \oplus X_2$.

Hence proceeding  by induction, we can construct a sequence
 $(u_n)_{n \in \N}$ of unit vectors in $\ell_2$  satisfying

\begin{itemize}

\item  $Au_{2n-1}=u_{2n}$ and $Au_{2n}=-u_{2n-1}$, $\forall n \in \N$;

 \item $e_n \in \sum_{k=1}^n  X_k$, $\forall n \in \N$;

 \item $u_i \perp u_j$, $\forall i, j \in \N$, $i \neq j$.
\end{itemize}

Thus the proof of proposition is complete. \pff
 
\

Note that from this proposition we may deduce that $\ell_2$ also admits a unique complex structure up to isomorphism, a well-known fact for which there does not seem to be a written proof in the literature. Indeed if $A$ is an operator on $\ell_2$ satisfying $A^2=-Id$, let $< , >'$ be the scalar product defined on $\ell_2$ by
$<x,y>'=<x,y>+<Ax,Ay>$ and $\norm{.}'$ be the associated norm.
Then $(\ell_2,\norm{.}')$ is a Hilbert space for which $A$ is an isometry, and therefore $(\ell_2^A,\norm{.}')$
is isometric to the canonical complex structure on $\ell_2$, hence $(\ell_2^A,\norm{.})$ is isomorphic to that complex structure.

\

In a first version of this paper we mentioned as an open question whether the
spaces $c_0$ and $\ell_p, p \neq 2$ admitted a unique complex structure up to
isomorphism. N.J. Kalton then indicated to us a nice and simple proof that this
is indeed the case. We reproduce this proof here with his authorization.

Recall that a Banach space $X$ is {\em primary} if $X \simeq Y$ or $X \simeq Z$
whenever $X=Y \oplus Z$.

\begin{lemm}\label{kaltonl} Let $X$ be a real Banach space, $A, B$ be  operators on $X$ such
  that $A^2=B^2=-Id$. Assume that $X^A$ is isomorphic to its conjugate and
  primary, and that $A$ and $B$ commute. Then $X^A$ and $X^B$ are isomorphic.
\end{lemm}

\pf It is easily checked that $P=\frac{1}{2}(Id+AB)$ and $Q=\frac{1}{2}(Id-AB)$
are projections on $X$ which commute with $A$ and $B$, and such that $Ax=-Bx$
for any $x \in PX$,  and
$Ax=Bx$ for any $x \in QX$. Let $X=Y \oplus Z$ be
the associated decomposition. Then $X^A \simeq Y^A \oplus Z^A$, and
$$X^B \simeq Y^B \oplus Z^B \simeq Y^{-A} \oplus Z^A.$$
Since $X^A$ is primary, we have either $Y^A \simeq X^A$ or $Z^A \simeq X^A$.
In the first case, and since $X^A$ is isomorphic to its conjugate, we deduce 
$$X^B \simeq Y^A \oplus Z^A \simeq X^A.$$
In the second case,
$$X^B \simeq Y^{-A} \oplus X^A \simeq Y^{-A} \oplus Z^{-A} \simeq X^{-A}
\simeq X^A.$$
\pff

\begin{theo}\label{kalton}(N.J.  Kalton) Let $X$ be a real Banach space  and
  assume that the complexification $X \oplus_{\C} X$ of $X$ is primary. Then $X$
  admits no more than one complex structure up to isomorphism. \end{theo}

\pf Let $J$ be the operator associated to the canonical complex structure on
$X \oplus X$, i.e. $J$ is defined by $J(x,y)=(-y,x)$. Assume $X$ admits a
complex structure and let $A$ be any
operator 
on $X$ such that $A^2=-Id$. Let $A \oplus A$ denote the corresponding operator
on $X \oplus X$. It is immediate that $A \oplus A$ and $J$ commute, therefore
 by Lemma \ref{kaltonl},
$$X^A \oplus X^A \simeq (X \oplus X)^{A \oplus A} \simeq (X \oplus X)^J.$$
Since the space $(X \oplus X)^J$ is primary we deduce that
$$X^A \simeq (X \oplus X)^J,$$
which concludes the proof. \pff

\begin{coro} The  spaces $c_0$, $C([0,1])$,  $\ell_p, 1 \leq p \leq +\infty$,
$L_p, 1 <p < +\infty$  admit a unique complex structure up to isomorphism.
\end{coro}

\subsection{Isometric complex structures up to renorming}

For other spaces than $\ell_2$ with its usual norm, the theory of complex structures up to isometry is quite different from the theory up to isomorphism.
For example, according to Jarosz in \cite{J}, any Banach space admits an equivalent norm for which the only isometries are $Id$ and $-Id$, and therefore does not admit complex structure in the isometric sense for that norm. We shall now use the methods of Jarosz to prove:

\begin{theo}\label{renormingcomplexstructure} Any real Banach space of dimension at least $4$ and which admits a complex structure up to isomorphism has an equivalent norm for which it admits exactly two complex structures up to isometry.

Any real Banach space which is isomorphic to a cartesian square has an equivalent norm for which it admits a unique complex structure up to isometry.
\end{theo}

Therefore $\ell_2$ with its usual norm is far from being the only Banach space with unique complex structure up to isometry. Actually Theorem \ref{renormingcomplexstructure} shows that all classical spaces may be renormed to have no, a unique, or exactly two complex structures up to isometry.
Note that the space $X(\C)$ defined in \cite{F}, which admits exactly two complex structures up to isomorphism, which are conjugate, cannot be renormed to admit a unique complex structure up to isometry. For if
$T^2=-Id$ for $T$ an isometry in some equivalent norm $\||.\||$, then
$X^T$ and $X^{-T}$ are complex structures on $X$ in the isomorphic sense, hence non-isomorphic and in particular non $\||.\||$-isometric.
Therefore the second part of Theorem \ref{renormingcomplexstructure} cannot be improved to all Banach spaces admitting a complex structure up to isomorphism.

\

Let $X$ be a complex space. We shall say that a real isometry $T$ on $X$
is {\em trivial} if $T=\lambda Id$, for some $\lambda$ in the complex unit circle.
We say that the complex linear structure on $X$ is {\em canonical} if $X$ is equal to a real cartesian sum $Y \oplus Y$ such that the multiplication by $i$ on $X$ is defined by $i(y,z)=(-z,y)$, for $y,z \in Y$ (i.e. $X$ is canonically isomorphic to the complexification of $Y$) and the conjugation map $c$  (i.e $c$ is defined by $c(y,z)=(y,-z)$) is an isometry on $X$. A real isometry $T$ on $X$ is then said to be a {\em conjugation isometry} if it is of the form
$T=\lambda c$, where $\lambda$ is in the unit complex circle.

\

The rest of this section is devoted to proving  1) that any complex Banach space of dimension at least $2$ may be renormed to admit only trivial real isometries, 2) that any complex Banach space real isomorphic to a cartesian square, and whose complex structure is canonical associated to that decomposition as a square, may be renormed to admit only trivial and conjugation real isometries, Corollary \ref{maintheorem}. Theorem \ref{renormingcomplexstructure} follows immediately
from Corollary \ref{maintheorem}. Indeed in case 1) the only isometries of square $-Id$ are $i Id$ and $-i Id$. Furthermore since the group of isometries commutes, $i Id$ and $-i Id$ are not conjugate in that group, so the associated complex structures are not isometric. There are therefore exactly two complex structures up to isometry,  which are conjugate.
In case 2), since $T^2=|\lambda|^2 Id$ whenever $T=\lambda c$, the isometries $i Id$ and $-i Id$ are also the unique isometries of square $-Id$. Since $-i Id=c(i Id)c=c^{-1}(i Id)c$, they are isometrically conjugate and their associated complex structures are $\C$-linearly isometric. Therefore there is a unique complex structure up to isometry in that case. 

\

Our proof consists in extending the methods of Jarosz concerning $\C$-linear isometries on complex spaces to the study of $\R$-linear isometries on complex spaces.

We first note that any equivalent norm on $\C$ is a multiple of the modulus; therefore real isometries on $\C$ are either trivial or conjugate in any equivalent norm, and $\C$ cannot be renormed to admit only trivial real isometries. We shall need a direct proof that the case of $\C^2$
is already different:

\begin{lemm}\label{C2} There exists a complex norm on $\C^2$ for which $\C^2$ only admits trivial real isometries.
\end{lemm}

\pf We fix $\lambda_0=0$ and $\lambda_k$, $1 \leq k \leq 4$, satisfying:

\begin{itemize}
\item i) $|\lambda_k|=1, \forall 1 \leq k \leq 4$,  

\item ii) $Re(\lambda_k)>0, \forall 1 \leq k \leq 4$,

\item iii) $\lambda_j \lambda_k \neq \lambda_l \lambda_m$ whenever 
$j,k,l,m \in \{1,2,3,4\}$ and $\{j,k\} \neq \{l,m\}$.
\end{itemize}

We define a norm $\norm{.}$ on $\C^2$ by the formula
$$\norm{(x,y)}=max \{|x|, max_{1 \leq k \leq 4} |x-\lambda_k y|\}=
max_{0 \leq k \leq 4}|x-\lambda_k y|,$$
and shall prove that any real isometry on $\C^2$ for that norm is trivial.

\

For $0 \leq k \leq 4$, let $$A_k=\{(x,y): |x-\lambda_k y|>max_{0 \leq j \leq 4, j \neq k}|x-\lambda_j y|\},$$ and let $$A=\cup_{0 \leq k \leq 4} A_k.$$
We observe that if $(x,y) \in A_k$, then $\norm{.}$ is constant in a neighborhood $(x,y)+V_k$ of $(x,y)$ in $(x,y)+H_k$, for some $H_k$ a linear subspace of $\R$-dimension $2$ (take $H_k=\{(\lambda_k h,h), h \in \C\}$); and note that $H_k$ is unique with this property.
On the other hand when $(x,y) \notin A$, let $j \neq k$ be such that
$\norm{(x,y)}=|x-\lambda_j y|=|x-\lambda_k y|$. If $\norm{.}$ is constant on a neighborhood $(x,y)+V$ of $(x,y)$ in $(x,y)+H$, for some $H$ a linear subspace of $\R$-dimension $2$, then $|x+h-\lambda_k(y+h')|$ is maximal for $(h,h')=0$ on $V$, from which we deduce that $h-\lambda_k h'=0$ for $(h,h') \in V$, so $H=H_k$. But then by the same reasoning $H=H_j$, a contradiction.
Finally we have proved that a point $(x,y)$ belongs to $A$ if and only if
 $\norm{.}$ is constant in a neighborhood of $(x,y)$ in $(x,y)+H$, for some $H$ a linear subspace of $\R$-dimension $2$; and so $A$ is defined by $\R$-linear and metric properties.

Let $T$ be an $\R$-linear isometry on $(\C^2,\norm{.})$. Therefore $T$ preserves $A$. Let
$(x,y) \in A_k, 0 \leq k \leq 4$ and let $0 \leq l \leq 4$ be such that
$T(x,y) \in A_l$. Since
$$\norm{T(x,y)}=\norm{(x,y)}=\norm{(x,y)+(h,h')}=\norm{T(x,y)+T(h,h')},$$
for $(h,h') \in V_k$, it follows that $T(H_k)=H_l$. We deduce easily that  $T(A_k) \subset A_l$. So $T(A_k)=A_l$ by symmetry, and finally there is a permutation $\sigma$ on $\{0,1,2,3,4\}$ such that $T(A_k)=A_{\sigma(k)}$ for all $0 \leq k \leq 4$. 



The isometry $T$ is given by a formula of the form
$$T(x,y)=(Ax+B\overline{x}+Cy+D\overline{y},ax+b\overline{x}+cy+d\overline{y}),$$
where $A,B,C,D,a,b,c,d$ are complex numbers.
It is easy to check that for any $0 \leq \theta \leq 2\pi$, 
$(e^{i\theta},0) \in \cap_{0 \leq k \leq 4}\overline{A_k}$ (for $k=0$ this requires condition ii)). By our computation of $T(A_k), 0 \leq k \leq 4$, we have therefore
$$T(e^{i\theta},0)=(Ae^{i\theta}+Be^{-i\theta}, ae^{i\theta}+be^{-i\theta}) \in \cap_{0 \leq k \leq 4}\overline{A_k},$$ and we deduce
$$1=\norm{T(e^{i\theta},0)}=|(A-\lambda_k a)e^{i\theta}+(B-\lambda_k b)e^{-i\theta}|, \forall 0 \leq \theta \leq 2\pi, \forall 0 \leq k \leq 4,.$$
We deduce easily that either $A=a=0$ in which case $|B-\lambda_k b|=|B|=1$ for all $1 \leq k \leq 4$, so also $b=0$; or that $B=b=0$ and similarly $|A|=1$, $a=0$.

Likewise for any $0 \leq \theta \leq 2\pi$, 
$(0,e^{i\theta}) \in \cap_{1 \leq k \leq 4}\overline{A_k}$.
 Then
$T(0,e^{i\theta})=(Ce^{i\theta}+De^{-i\theta}, ce^{i\theta}+de^{-i\theta}) \in \cap_{0 \leq k \leq 4, k \neq \sigma(0)}\overline{A_k}$, so we deduce
$$1=|(C-\lambda_k c)e^{i\theta}+(D-\lambda_k d)e^{-i\theta}|,\forall 0 \leq \theta \leq 2\pi, \forall 0 \leq k \leq 4, k \neq \sigma(0).$$
So either $C=c=0$ in which case $|D-\lambda_k d|=1$ for all $1 \leq k \leq 4$, $k \neq \sigma(0)$, from which it follows easily that $d=0$ or $D=0$; or $D=d=0$ and ($C=0$ or $c=0$).

\

Summing up we have obtained that $T$ is given either by 1) $T(x,y)=(Ax,cy)$, 2) $T(x,y)=(Ax,d\overline{y})$, 3) $T(x,y)=(B\overline{x},cy)$, or 4) $T(x,y)=(B\overline{x},d\overline{y})$. It remains to prove that only 1) is possible, with $A=c$. We may assume $A=B=1$, and we have $|c|=|d|=1$.

Fixing $1 \leq j \leq 4$ and $0 \leq \theta \leq 2\pi$, we observe that
$$2=\norm{(e^{i\theta},-\overline{\lambda_j} e^{i\theta})}.$$
If $T$ satisfies 2), then we deduce
$$2=\norm{(e^{i\theta},-d\lambda_j e^{-i\theta})}=
max(1, max_{1 \leq k \leq 4}|e^{i\theta}+\lambda_k\lambda_j de^{-i\theta}|).$$
 So $2=max_{1 \leq k \leq 4}|e^{i2\theta}+\lambda_k\lambda_j d|$, but obviously this is only possible for a finite number of values of $\theta$, so we get a contradiction. A similar reasoning holds to exclude the case 3).

 If $T$ satisfies 4), then
$$2=\norm{(e^{-i\theta},-d\lambda_j e^{-i\theta})}=
max(1, max_{1 \leq k \leq 4}|e^{-i\theta}+\lambda_k\lambda_j de^{-i\theta}|).$$
 So $2=max_{1 \leq k \leq 4}|1+\lambda_k\lambda_j d|$. We deduce that $\forall 1 \leq j \leq 4$, $\exists 1 \leq k \leq 4: \lambda_j \lambda_k = 1/d$; but this contradicts condition iii) on the $\lambda_k$'s.

So $T$ satisfies 1) and it remains to prove that $c=1$. We have that
 $$2=\norm{(e^{i\theta},-c\overline{\lambda_j} e^{i\theta})}=
max(1, max_{1 \leq k \leq 4}|e^{i\theta}+\lambda_k\overline{\lambda_j} ce^{i\theta}|).$$
So $2=max_{1 \leq k \leq 4}|1+\lambda_k\overline{\lambda_j} c|$. 
We deduce that $\forall 1 \leq j \leq 4$, $\exists 1 \leq k \leq 4: \lambda_j= c \lambda_k$. But then $c=1$, otherwise there exist $2 \leq k \leq 4$ and $1 \leq k' \leq 4$ such that $\lambda_1/\lambda_k=c=\lambda_k/\lambda_{k'}$, contradicting condition iii).
\pff

\begin{prop}\label{double} Let $\Gamma$ be a nonempty set and $E$ a complex Banach space such that $c_0(\Gamma,\C) \subset E \subset \ell_{\infty}(\Gamma,\C)$. Then 

1) if $|\Gamma| \geq 2$ then there is a norm $\||.\||_1$ on $E$, equivalent with the original sup norm $\|.\|$ of $E$ and such that an $\R$-linear map $T$ on $E$ is both a $\|.\|$ and $\||.\||_1$ isometry if and only if $T$ is trivial;

2) if $E=X \oplus X$, where $X$ is a real Banach space such that
$c_0(\Gamma,\R) \subset X \subset \ell_{\infty}(\Gamma,\R)$, and $E$ is equipped with the norm induced by the canonical inclusion of $E$ in
$l_{\infty}(\Gamma,\R) \oplus \ell_{\infty}(\Gamma,\R)$ identified with $l_{\infty}(\Gamma,\C)$, then
there is a norm $\||.\||_2$ on $E$, equivalent with the original norm $\|.\|$ of $E$ and such that an $\R$-linear map $T$ on $E$ is both a $\|.\|$ and $\||.\||_2$ isometry if and only if $T$ is a trivial or conjugate isometry.
\end{prop}

\pf We adapt a proof of Jarosz, \cite{J} Proposition 1. We first note that vectors $x,y$ in $\ell_{\infty}(\Gamma,\C)$ do not have disjoint supports if and only if 
$$(1) \exists
z \in \ell_{\infty}, \exists \epsilon=\pm 1: \norm{z} \leq 1, \norm{x+z} \leq 1, \norm{\epsilon y+z} \leq 1,  \norm{x+\epsilon y+z}>1.$$
The if part follows from the fact that whenever $x,y$ have disjoint supports on $\ell_{\infty}(\Gamma,\C)$, then for all $z \in \ell_{\infty}(\Gamma,\C)$, and for any $\epsilon=\pm 1$, 
$$\norm{x+\epsilon y +z} \leq max(\norm{z},\norm{x+z},\norm{\epsilon y +z}).$$
We prove the only if part.
Let $k \in supp(x) \cap supp(y)$. Up to replacing $y$ by $-y$, we may assume that the scalar product
$<x_k,y_k>=Re(x_k \overline{y_k})$ is positive, and note that $|x_k| \leq 1, |y_k| \leq 1$. Let $\lambda_1 \geq 0$ be such that $|x_k+\lambda_1 x_k|=1$, $\lambda_2 \geq 0$ be such that $|y_k+\lambda_2 x_k|=1$ and let
$\lambda=min(\lambda_1,\lambda_2)$. Then by easy geometrical considerations in $\R^2$ and the fact that $<x_k,y_k>$ is non-negative,
$|x_k+\lambda x_k| \leq 1$ and $|y_k+\lambda x_k| \leq 1$, but
$|x_k+y_k+\lambda x_k|>1$. It follows, letting $z=\lambda x_k e_k$, that $\norm{z} \leq 1$, $\norm{x+z} \leq 1, \norm{y+z} \leq 1$, but $\norm{x+y+z}>1$.

\

Now since (1) is preserved by $\R$-linear isometries, whenever $T$ is an $\R$-linear isometry for $\norm{.}$, we see that $T$ maps disjointly supported vectors to disjointly supported vectors. It follows easily that $\forall \gamma \in \Gamma$, $T$ maps $e_{\gamma}$ to some $\epsilon_{\gamma} e_{\pi(\gamma)}$, where $\pi:\Gamma \rightarrow \Gamma$ and 
$|\epsilon_{\gamma}|=1$. Since $Ti$ is also an $\R$-linear isometry, it also follows that $T$ maps $ie_{\gamma}$ to some $\epsilon_{\gamma}^{\prime} ie_{\pi'(\gamma)}$, where $\pi':\Gamma \rightarrow \Gamma$ and 
$|\epsilon^{\prime}_{\gamma}|=1$. Since for all $\theta \in \R$,
$$1=\norm{T(e^{i\theta}e_{\gamma})}=\norm{\cos\theta \epsilon_{\gamma} e_{\pi(\gamma)}+i\sin\theta  \epsilon_{\gamma}^{\prime} e_{\pi'(\gamma)}},$$
it follows easily that $\pi'(\gamma)=\pi(\gamma)$, that $\epsilon'_{\gamma}=\pm \epsilon_{\gamma}$, and therefore that $\pi$ is a bijection on $\Gamma$.

\

Observe now that to prove that an $\R$-linear isometry $T$ is equal to $\lambda Id$, it is enough to prove that $T(e_\gamma)=\lambda e_{\gamma}$ and $T(ie_\gamma)=i\lambda e_{\gamma}$ for all $\gamma \in \Gamma$. Indeed then for any $(a_\gamma)_\gamma \in E$, writing
$(b_\gamma)_\gamma=T((a_\gamma)_\gamma)$, we have for $r \in \C$ and fixed $\gamma$,
$$\norm{(a_\gamma)_\gamma - re_\gamma}=\norm{(b_\gamma)_\gamma - \lambda re_\gamma}.$$
When $|r|$ is large enough this means $|b_\gamma - \lambda r|=|a_\gamma -r|$ and therefore $b_\gamma=\lambda a_\gamma$.

For 2), a similar reasoning shows that it is enough to prove that  $T(\mu e_\gamma)=\lambda \overline{\mu} e_\gamma$ for all $\gamma \in \Gamma$ and $\mu \in \C$ to obtain that for any $(a_\gamma)_\gamma \in E$, $T((a_\gamma)_\gamma)=\lambda ((\overline{a_\gamma})_\gamma)$.

\

To prove 2), we may assume that $|\Gamma| \geq 2$. We fix a well-order $<$ on $\Gamma$ and define

$$\norm{x}_2=max(\norm{x}, \sup|2x(\gamma)+x(\beta)|, \gamma<\beta \in \Gamma).$$
Assume $T$ is an $\R$-linear isometry for $\norm{.}$ and $\norm{.}_2$.
If $\gamma<\gamma'$ but $\pi(\gamma)>\pi(\gamma')$ then 
$\norm{2e_\gamma+e_{\gamma'}}_2=5$ but
$\norm{T(2e_\gamma+e_{\gamma'})}_2=\norm{\epsilon_{\gamma'} e_{\pi(\gamma')} + 2\epsilon_\gamma e_{\pi(\gamma)}}_2 \leq 4$, a contradiction. So
 $\pi$ preserves order and is therefore equal to $Id_{\Gamma}$.
If $\epsilon_{\gamma} \neq \epsilon_{\gamma'}$ for $\gamma < \gamma'$ then $\norm{e_{\gamma}+e_{\gamma'}}_2=3$ but
$$\norm{T(e_{\gamma}+e_{\gamma'})}_2=\norm{\epsilon_{\gamma} e_{\gamma}+\epsilon_{\gamma'} e_{\gamma'}}_2 \leq max\{1,2,|2\epsilon_\gamma+\epsilon_{\gamma'}|\}<3.$$
Hence $\epsilon_{\gamma}$ is constant on $\Gamma$.

We have finally obtained that for some $\lambda=\pm \mu, |\lambda|=1$, and for all $\gamma \in \Gamma$, $T(e_\gamma)=\lambda e_\gamma$ and
$T(i e_\gamma)=\mu ie_\gamma$.
If $\lambda=\mu$ then we deduce that $T$ is a trivial isometry and if $\lambda=-\mu$ then $T$ is a conjugate isometry.

\

To prove 1), we fix some $\gamma_0<\gamma_1$ and consider the norm defined by
$$\norm{x}_1=max(\norm{x}_2, |\frac{3}{2}x(\gamma_0)+ix(\gamma_1)|).$$
It is easy to check that we may repeat the reasoning used for $\norm{.}_2$ to obtain that if $T$ is an $\R$-linear isometry for $\norm{.}$ and
$\norm{.}_1$, then for some $\lambda=\pm \mu$, $|\lambda|=1$, and for all $\gamma \in \Gamma$, $T(e_\gamma)=\lambda e_\gamma$ and
$T(i e_\gamma)=\mu i e_\gamma$. Furthermore, since
$$\norm{e_{\gamma_0}+ie_{\gamma_1}}_1=max(1,2,\sqrt{5},\frac{1}{2})=\sqrt{5},$$
but 
$$\norm{e_{\gamma_0}-ie_{\gamma_1}}_1=max(1,2,\sqrt{5},\frac{5}{2})=
\frac{5}{2} \neq \sqrt{5},$$
we deduce that $\lambda=\mu$ and therefore $T$ is a trivial isometry. 
\pff

\

Observe that if $|\Gamma|=1$, then $E=\C$, but it is clearly not possible to renorm $\C$ to admit only trivial real isometries. So the condition that $|\Gamma| \geq 2$ in Proposition \ref{double} is necessary. 

The next proposition is a version of Proposition 3 from \cite{J} for real isometries on complex spaces. A great part of its proof is identical to the proof of \cite{J} Proposition 3. However, we were not able to prove the affirmation that "evidently $Tx_0=x_0$" in the original proof of Jarosz. Our version is  therefore a weaker version of \cite{J} Proposition 3 in the case of real isometries on complex spaces, in which we added the hypotheses 1) or 2); this weaker version is enough for the proof of Theorem \ref{qwerty}. Note that it is clear by the same reasoning that Proposition 3 from \cite{J} is valid if one adds in the hypothesis that for any $T$ which is an isometry for $\norm{.}$ and $p(.)$, $Tx_0$ and $x_0$ are linearly dependent, and that this is enough to deduce \cite{J} Theorem 1.

\begin{prop}\label{2cas} Let $(X,\|.\|)$ be a complex Banach space, $x_0$ a non-zero element of $X$, $p(.)$ a continuous norm on $(X,\|.\|)$, $G_1$ the group of all real isometries of $(X,\|.\|)$ and $G_2$ the group of all real isometries of $(X,p(.))$. Assume

1) for any $T \in G_1 \cap G_2$, there exists $\lambda \in C$ such that
$Tx_0=\lambda x_0$ and $Tix_0=\lambda ix_0$, or

2) the linear structure on $X$ is canonical, and for any $T \in G_1 \cap G_2$, there exists $\lambda \in C$ such that
$Tx_0=\lambda x_0$ and $Tix_0=\pm \lambda ix_0$,

Then
there is a norm $\norm{.}_w$ on $Y=X \oplus \C$ such that
$\norm{.}_w$ and $\norm{.}$ coincide on $X$ and
the group $G$ of real isometries of $(Y,\norm{.}_w)$ is isomorphic to $G_1 \cap G_2$.
\end{prop}

\pf The beginning of the proof is the same for 1) and 2), then we shall differentiate the proof at the end. In 1), the isomorphism $\alpha_1$ from 
$G$ onto $G_1 \cap G_2$ is the restriction map $\alpha_1(T)=T_{|X}$ and its inverse is given by $\alpha_1^{-1}(T)=T \oplus \lambda Id_{\C}$, if $\lambda$ is such that $Tx_0=\lambda x_0$ and $Tix_0=\lambda ix_0$.
In 2), the isomorphism $\alpha_2$ from 
$G$ onto $G_1 \cap G_2$ is also the restriction map and the inverse is given by
$\alpha_2^{-1}(T)=T \oplus \lambda Id_{\C}$, if $\lambda$ is such that $Tx_0=\lambda x_0$ and $Tix_0=\lambda ix_0$, or
$\alpha_2^{-1}(T)=T \oplus \lambda c_{\C}$, if $\lambda$ is such that $Tx_0=\lambda x_0$ and $Tix_0=-\lambda ix_0$ (here $c_{\C}$ is the conjugation map on $\C$).

By replacing $p(.)$ by $p(.)+\norm{.}$ and multiplying by an appropriate number, we may assume that $p$ and $\norm{.}$ are equivalent, that $1000\norm{.} \leq p(.)$ and that $\norm{x_0} \leq 0.1$.
Let 
$$A=\{(x,\alpha) \in X \oplus \C=Y: max\{\norm{x},|\alpha|\} \leq 1\},$$
$$C=\{(x+x_0,2): p(x) \leq 1\},$$
and let $\norm{.}_W$ be the norm whose unit ball $W$ is the closed balanced convex set generated by $A \cup C$. Observe that 
$\norm{(x,\alpha)}_W=\norm{x}$ whenever $|\alpha| \leq \norm{x}$. Therefore the norm $\norm{.}_W$ coincides with the original one on $X$. It is also evident that if $T: X \rightarrow X$ preserves both norm $\norm{.}$ and $p(.)$, $Tx_0=\lambda x_0$ and $Tix_0=\lambda ix_0$, $|\lambda|=1$ then
$T \oplus \lambda Id_{\C}$ is an isometry of $Y$.
In case 2) it also easy to check that if $T: X \rightarrow X$ preserves both norm $\norm{.}$ and $p(.)$, $Tx_0=\lambda x_0$ and $Tix_0=-\lambda ix_0$, $|\lambda|=1$ then
$T \oplus \lambda c_{\C}$ is an isometry of $Y$.

Assume now that $T: Y \rightarrow Y$ is a $\norm{.}_W$-isometry.
We first prove that $T$ maps $X$ onto $X$ and $T_{|X}$ preserves both $\norm{.}$ and $p(.)$.

We note that $C$ as well as all its rotations $\lambda C, |\lambda|=1$  are faces of $W$. 
We distinguish two types of points in $\delta W$:
A) points interior to a segment $I$ contained in $\delta W$, whose length (with respect to the $W$-norm) is at least $0.1$, and the limits of such points; B) all other points.

As these types are $\R$-linearly metrically defined, they are preserved by $T$.
On the other hand it is easy to see that the points of type A) cover all of $\delta W$ except the relative interiors of the faces $\lambda C$. Hence $T(x_0,2)$ belongs to some $\lambda C$ with 
$|\lambda|=1$.
Replacing $T$ by $\lambda^{-1}T$ we can assume that $T(x_0,2) \in C$ and since $T$ maps the face $C$ onto a face of $W$ we have $TC=C$. To prove that $T$ maps $X$ onto $X$, let $x \in X$ with $p(x) \leq 1$. We have 
$$T(x,0)=T((x+x_0,2)-(x_0,2))=T(x+x_0,2)-T(x_0,2) \in C - C \subset X,$$
and as $\{x: p(x) \leq 1\}$ contains a ball in $X$ this is true for all $x \in X$, i.e. $TX \subset X$; by symmetry $TX=X$. Because the $\norm{.}_W$ norm agrees with $\norm{.}$ on $X$, it follows that $T_{|X}$ is a 
$\norm{.}$-isometry.
Since $TC=C$ the function $T_{|X}$ maps $(x_0,2)$ onto itself and the unit ball for $p(.)$ onto itself. Therefore $T_{|X}$ preserves the norm $p(.)$ as well.

We prove a fact that will allow us to conclude the proof in cases 1) and 2): if $T$ is an isometry on $Y$ for $\norm{.}_W$, and $\lambda$ is such that $T(x_0,0)=\lambda (x_0,0)$, then $T(0,1)=(0,\lambda)$. To check this
we may assume $\lambda=1$. By the above we have that $T(x_0,2)=\mu(x_0,2)$ for some $\mu, |\mu|=1$. Therefore by an easy computation, for all $k \in \R$,
$$T(kx_0,1)=((k+\frac{\mu-1}{2})x_0,\mu).$$
If $\mu \neq 1$ then pick $k=\pm \frac{1}{\norm{x_0}}$ so that
$|k+\frac{\mu-1}{2}|>\frac{1}{\norm{x_0}}$. Then $(kx_0,1)$ belongs to $W$ while $T(kx_0,1)$ does not, a contradiction. Therefore $\mu=1$, i.e.
$T(x_0,2)=(x_0,2)$ and so $T(0,1)=(0,1)$. Therefore the fact is proved.

In case 1), considering $T$ an isometry on $Y$ for $\norm{.}_W$, we know that $T_{|X}$ belongs to $G_1 \cap G_2$, therefore there exists $\lambda$ such that $Tx_0=\lambda x_0$ so by the fact, $T(0,1)=(0,\lambda)$. We also have that $Tix_0=\lambda ix_0$, so by the fact for $Ti$,
$T(0,i)=(0,\lambda i)$. Finally $T(0,z)=(0,\lambda z)$ for all $z \in \C$.

In case 2), again an isometry $T$ on $Y$ for $\norm{.}_W$ is such  that $T_{|X}$ belongs to $G_1 \cap G_2$. If $Tx_0=\lambda x_0$ and $Tix_0=\lambda ix_0$ then as above $T(0,z)=(0,\lambda z)$ for all $z \in \C$.
If $Tx_0=\lambda x_0$ and $Tix_0=-\lambda ix_0$ then by the fact applied to $T$ and $Ti$, $T(0,z)=(0,\lambda \overline{z})$ for all $z \in \C$.
\pff

\

By the form of the isomorphisms between $G$ and $G_1 \cap G_2$ given at the beginning of the proof of Proposition \ref{2cas} in cases 1) and 2), we deduce immediately:

\begin{coro}\label{corol} Let $(X,\|.\|)$ be a complex Banach space, $p(.)$ a continuous norm on $(X,\|.\|)$. Then
there is a norm $\norm{.}_w$ on $Y=X \oplus \C$ such that
$\norm{.}_w$ and $\norm{.}$ coincide on $X$ and such that:

1) If every real isometry on $X$ for $\norm{.}$ and for $p(.)$ is trivial, then every real isometry on $Y$ for $\norm{.}_w$ is trivial.

2) If the linear structure on $X$ is canonical, and every real isometry on $X$ for $\norm{.}$ and for $p(.)$ is a trivial or a conjugacy isometry, then every real isometry on $Y$ for $\norm{.}_w$ is a trivial or a conjugacy isometry. 
\end{coro}

\

The following fact, due to Pli\v cko \cite{P}, was cited and used in \cite{J}.

\begin{prop}\label{plicko}(Pli\v cko \cite{P}) For any Banach space $X$ there is a set $\Gamma$ and a continuous, linear injective map $J$ from $X$ into $\ell_{\infty}(\Gamma)$ such that the closure of $J(X)$ contains $c_0(\Gamma)$.
\end{prop}

\begin{theo}\label{qwerty} For any complex Banach space $X$ of dimension at least $1$ (resp. any complex Banach space $X$ with canonical linear structure), there is a Banach space $Y$ with $X \subset Y$ and $\dim Y/X=1$ such that $Y$ has only trivial real isometries (resp. trivial and conjugation real isometries).
\end{theo}

\pf It imitates the proof of \cite{J} Theorem 1. Let $Y=X \oplus \C$. If $\dim X=0$ then the result is trivial. If $\dim X=1$ in the trivial isometries case, then the proof holds from Lemma \ref{C2}, since we may assume that the norm on $X \simeq \C$ is the modulus, and the norm
$\norm{.}$ on $\C^2$ from Lemma \ref{C2} satisfies
$\norm{(x,0)}=|x|$, for all $x \in \C$. 

In other cases,
let $J: X \rightarrow \ell_{\infty}(\Gamma,\C)$ be an injective map given by Proposition \ref{plicko}; in the case when $X$ has canonical linear structure, $X=Z \oplus Z$, then define an injective map
$j: Z \rightarrow \ell_{\infty}(\Gamma,\R)$ and $J=j \oplus j:
X \rightarrow \ell_{\infty}(\Gamma,\C)$.
Let $E:=\overline{J(X)} \subset \ell_{\infty}(\Gamma,\C)$. We prove that 
there is a continuous norm $\tilde{p}$ on $E$ such that $(E,\tilde{p})$ has only trivial (resp. trivial and conjugacy) isometries. If $\dim X=1$ in the trivial or conjugacy isometries case, then this is obvious. If $\dim X=2$ in the trivial isometries case, this holds by Lemma \ref{C2}. 
If now  $\dim X \geq 2$ in the trivial or conjugacy isometries case, or  $\dim X \geq 3$ in the trivial isometries case, then fix $\gamma \in \Gamma$.
Let $\||.\||$ be the norm $\norm{.}_1$ (resp. $\norm{.}_2$) on 
$\{e \in E: e(\gamma)=0\} \subset \ell_{\infty}(\Gamma \setminus \{\gamma\},\C)$ given by Proposition \ref{double}. We then have
$$E \simeq \{e \in E: e(\gamma)=0\} \oplus_{\infty} \C,$$
so by Corollary \ref{corol}, there is a continuous norm $\tilde{p}$ on $E$ such that $(E,\tilde{p})$ has only trivial (resp. trivial and conjugacy) isometries. We define a continuous norm $p$ on $X$ by
$$p(x)=\tilde{p}(Jx), \ \ \ x \in X.$$
Evidently $(JX,\tilde{p})$ and so $(X,p)$ have only trivial (resp. trivial and conjugacy) isometries. Hence, again by Corollary \ref{corol}, there is a norm on $Y=X \oplus \C$, with only trivial (resp. trivial and conjugacy) isometries, which coincides with $\norm{.}$ on $X$.\pff

\begin{coro}\label{maintheorem} For any complex Banach space $X$ of dimension at least $2$ (resp. any complex Banach space $X$ with canonical linear structure), there is an equivalent norm on $X$ for which $X$ has only trivial (resp. trivial and conjugation) real isometries.
\end{coro}

\section{Isometries on real HI spaces}

It may be interesting to conclude 
this article by noting that 
isometries on the real HI space of Gowers and Maurey, or more generally, on spaces with the $\lambda Id +S$ property, have specific properties under any equivalent norm. This was obtained in \cite{RR} in the complex case.

 Recall that the complexification $\hat X$ of a real Banach space $X$ (see, for example, \cite{LT}, page 81) is defined as the space 
 $\hat X= \{x+iy: x, y \in X \},$ which is the space $X \oplus X$ with the canonical complex structure associated to $J$ defined on $X \oplus X$ by $J(x,y)=(-y,x)$.
 Let $A, B \in L(X)$. Then 
 $$(A+iB)(x+iy):=Ax-By+(Ay+Bx)$$
 defines an operator $A+iB \in L(\hat X)$ that satisfies 
 $\max \{ \|A \|, \|B \| \} \leq \|A+iB \| \leq 2^{1/2} (\|A \|+ \|B\|).$
 Conversely, given $T \in L(\hat X)$, if we put $T(x+i0):=Ax+iBx$, then we obtain $A, B \in L(X)$ such that $T=A+iB$.
 We write $\hat{T}=T+i0$ for $T \in {\cal L}(X)$.

\

Let $T \in L(X)$. We recall that the group $(e^{tT})_{t \in \R}$ has growth order $k \in \N$ if
 $\|e^{tT} \|=\sigma(|t|^k) \ \ \hbox{as}  \ \  |t| \to + \infty.$
 We also recall that an invertible operator $T \in L(X)$ is polynomially bounded of order $k \in \N$ if
 $\|T^n \|=\sigma(n^k) \ \ \hbox{as}  \ \  |n| \to + \infty.$
 In \cite{RR}, Theorem 3.2, it is proved that:

\begin{prop}\cite{RR} Let $X$ be a complex Banach space and $T \in L(X)$ such that there exists $\lambda \in \C$ with $T- \lambda I \in S(X)$ and the group     $(e^{tT})_{t \in \R}$     has growth order $k \in \N$. Then $(T- \lambda I)^k$ is a compact operator.
 \end{prop}

 The result in \cite{RR} is stated for complex HI spaces but the proof only uses the fact that complex HI spaces satisfy the $\lambda Id +S$-property.
 So by using this proposition instead of \cite{RR} Theorem 3.2, we can prove in similar way to \cite{RR} Theorem 3.5 the following result:

 \begin{prop}\label{RRR} Suppose that $X$ is a complex Banach space with the $\lambda Id +S$ property and $T \in L(X)$ is an  invertible operator, polynomially bounded of order $k \in \N$. Let $\lambda \in \C$  such that  $T- \lambda I \in S(X)$.  Then $(T- \lambda I)^k$ is a compact operator.
 \end{prop}
 
 We deduce:

 \begin{prop}\label{mnbvc} Suppose that  $X$ is a real Banach space with the $\lambda Id +S$-property and $T \in L(X)$  is an isometry. Then $T$ is of the form $\pm Id+K$,  $K$ compact.
 \end{prop}

 \pf Let $T$ be an isometry on $X$, and $a \in \R$, $S$ strictly singular be such that $T=a Id+S$. Clearly $a=\pm 1$. Let $\hat{X}$ be the complexification of $X$. Using \cite{Gon2} Proposition 2.6, it is easy to check that $\hat{X}$ has the $\lambda Id+S$-property. Consider $\hat{T}=T+i.0 \in L(\hat{X})$. Notice that $\hat{T}$ is an isomorphism from   $\hat{X}$ onto $\hat{X}$. Moreover,  ${\hat{T}}^{n}=T^{n}+i.0, \ \ \forall n \in \N$ and thus $\|{\hat{T}}^{n} \|$ is bounded. In particular, $\hat{T}$ is polynomially bounded of order 1.

Now notice that $\hat{T} - a \hat{Id}=(T-aId)+i.0$. Thus by \cite{Gon2} Proposition 2.6, $\hat{T} - a \hat{Id} \in S(\hat{X})$. Therefore according to  Proposition \ref{RRR}, $\hat{T} - \lambda Id$ is a compact operator. So by \cite{Gon2} Proposition 2.4,  $T-aId$ is also a compact operator.

\

\begin{ques} Let $X$ be a real H.I. Banach space such that
every operator is of the form $\lambda Id +\mu J +S$, where $J^2=-Id$. Does it follow that every isometry is of the form $\lambda Id +\mu J +K$, $K$ compact?  
\end{ques}

In this direction, it is natural to ask whether the complexification of a real HI space is always HI. By the proof of \cite{F} Proposition 35 this is always the case when every operator on a subspace $Y$ of $X$ is of the form $\lambda i_{YX}+s$, $\lambda \in \R$, $s$ strictly singular. Observe that by \cite{Gon2} Proposition 3.16,  if a real Banach space $X$ is such that ${\cal L}(X)/{\cal S}(X)$ is isomorphic to $\C$ or $\HH$, then the complexification $\hat{X}$ of $X$ is decomposable.
Therefore if $\hat{X}$ is HI for some real space $X$, then every subspace of $X$ must have the $\lambda Id+S$-property.

\section{Appendix}
We give the proof of two lemmas used in Section 2. They are inspired by \cite{DGZ} Theorem 7.4 page 72 and by the properties of Day's norm on $c_0$, \cite{DGZ} page 69.

\begin{lemm}\label{appendix} Let $Y$ be a Banach space with an LUR norm, let $1 \leq p <+\infty$, and let $X=l_p(Y)$. Then there exists an equivalent LUR norm on $X$ for which any map $T$ defined on $X$
by $T((y_n)_{n \in \N})=(\epsilon_n y_{\sigma(n)})_{n \in \N},$ where
$\epsilon_n=\pm 1$ for all $n \in \N$ and $\sigma$ is a permutation on $\N$, is an isometry.
\end{lemm}

\pf Fix an equivalent LUR norm $\norm{.}$ on $Y$, and let $\norm{.}=\norm{.}_p$ be the corresponding $l_p$-norm on $X$, when $p>1$. When $p=1$, let
$\norm{.}_1$ denote the corresponding $l_1$-norm,
$\norm{.}_2$ denote the corresponding $l_2$-norm (via the canonical "identity" map from $l_1$ into $l_2$), and let $\norm{.}$ be the equivalent norm defined on $X$ by
$$\norm{x}^2=\norm{x}_1^2+\norm{x}_2^2.$$
To prove that $\norm{.}$ is LUR let $x=(y_k)_k \in X$ and $x_n=(y_{n,k})_k \in X$ with $\lim_n \norm{x_n}=\norm{x}$ and 
$\lim_n \norm{x+x_n}=2\norm{x}$. We need to prove that 
$\lim_n x_n=x$.

We first assume that $p=1$. We have that $$\lim_n 2\norm{x}^2+2\norm{x_n}^2-\norm{x+x_n}^2=0. \eqno (1)$$
Using \cite{DGZ} Fact 2.3 p 45, (1) implies
$$\lim_n 2\norm{x}_1^2+2\norm{x_n}_1^2-\norm{x+x_n}_1^2=0 \eqno (2)$$
and
$$\lim_n 2\norm{x}_2^2+2\norm{x_n}_2^2-\norm{x+x_n}_2^2=0. \eqno (3)$$
By \cite{DGZ} Fact 2.3 again, (3) implies, for all $k \in \N$,
$$\lim_n 2\norm{y_k}^2+2\norm{y_{n,k}}^2-\norm{y_k+y_{n,k}}^2=0,$$
whence, since the norm on $Y$ is LUR, by \cite{DGZ} Proposition 1.2. p 42,
$$\lim_{n} y_{n,k}=y_k, \forall k \in \N, \eqno (4)$$
and from (2) we have, see \cite{DGZ} p 42,
$$\lim_n \norm{x_n}_1=\norm{x}_1. \eqno (5)$$

Now assume $p>1$. We have that $$\lim_n\norm{x_n}_p=\norm{x}_p \eqno (6)$$ which means that
$$\lim_n \sum_k \norm{y_{n,k}}^p= \sum_k \norm{y_k}^p. \eqno (7)$$
Let $|.|_p$ also denote the norm on $\ell_p$. Since
$$\norm{x_n+x}=|(\norm{y_{n,k}+y_k})_k|_p
\leq |(\norm{y_{n,k}}+\norm{y_k})_k|_p$$
$$\leq |(\norm{y_{n,k}})_k|_p+|(\norm{y_k})_k|_p
=\norm{x_n}+\norm{x}$$
and both $\norm{x_n+x}$ and $\norm{x_n}+\norm{x}$ converge to $2\norm{x}$,
we deduce that
$$\lim_n |(\norm{y_{n,k}}+\norm{y_k})_k|_p=2 |(\norm{y_k})_k|_p. \eqno (8)$$
Since $|.|_p$ is LUR on $\ell_p$, we deduce from (7) and (8) that
$\lim_n |(\norm{y_{n,k}}-\norm{y_k})_k|_p=0$, in particular
$$\forall k \in \N, \lim_n \norm{y_{n,k}}=\norm{y_k}.\eqno (9)$$
Since $\norm{x+x_n}$ converges to $2\norm{x}$ we also have
$$\lim_n \sum_k \norm{y_{n,k}+y_k}^p=2^p \sum_k \norm{y_k}^p. \eqno (10)$$
Fix $k_0 \in \N$ and $\epsilon>0$. We may find some $k_1>k_0$ such
that $$\sum_{k \geq k_1}\norm{y_k}^p<\epsilon. \eqno (11)$$
Therefore by (7), (9), and (11), for $n$ large enough,
$$\sum_{k \geq k_1}\norm{y_{n,k}}^p<2\epsilon. \eqno (12)$$
Using (9), (11) and (12), we deduce that for $n$ large enough,
$$\sum_k \norm{y_{n,k}+y_k}^p < 2^p \sum_{k \neq k_0, k<k_1} \norm{y_k}^p +\epsilon+2^p. 3\epsilon+\norm{y_{n,k_0}+y_{k_0}}^p, \eqno (13)$$
while by (10) and (11), for $n$ large enough,
$$\sum_k \norm{y_{n,k}+y_k}^p > 2^p \sum_{k \neq k_0, k<k_1} \norm{y_k}^p +2^p\norm{y_{k_0}}^p-2^p\epsilon-\epsilon. \eqno (14)$$
From (13) and (14) we deduce that for $n$ large enough,
$$2^p \norm{y_{k_0}}^p < (2+4.2^p)\epsilon+\norm{y_{n,k_0}+y_{k_0}}^p,$$
and we deduce, using also (9), that
$$\lim_n \norm{y_{n,k_0}+y_{k_0}}=2\norm{y_{k_0}}. \eqno (15)$$
From (9) and (15), and from the fact that the norm on $Y$ is LUR, it follows that
$$\forall k \in \N, \lim_n y_{n,k}=y_k. \eqno (16)$$

Going back to the general case, fix $\epsilon>0$ and let $k_1 \in \N$ be such that $\sum_{k \geq k_1}\norm{y_k}^p<\epsilon$, then
$$\norm{x-x_n}_p^p=\sum_{k<k_1} \norm{y_k-y_{n,k}}^p+\sum_{k \geq k_1}
\norm{y_k-y_{n,k}}^p$$
$$\leq \sum_{k <k_1} \norm{y_k-y_{n,k}}^p+2^p\sum_{k \geq k_1}
\norm{y_k}^p+2^p \sum_{k \geq k_1}
\norm{y_{n,k}}^p$$
$$=\sum_{k <k_1}\norm{y_k-y_{n,k}}^p+2^p(2\sum_{k \geq k_1}\norm{y_k}^p
+(\norm{x_n}_p^p-\norm{x}_p^p)+\sum_{k < k_1}(\norm{y_k}^p-\norm{y_{n,k}}^p)).$$
So by (4) and (5) when $p=1$, or by (6) and (16) when $p>1$, we obtain that
$\norm{x-x_n}_p^p<3.2^p\epsilon$ for $n$ large enough. \pff

\begin{lemm}\label{appendixbis}Let $Y$ be a Banach space with an LUR norm and let $X=c_0(Y)$. Then there exists an equivalent LUR norm on $X$ for which any map $T$ defined on $X$
by $T((y_n)_{n \in \N})=(\epsilon_n y_{\sigma(n)})_{n \in \N},$ where
$\epsilon_n=\pm 1$ for all $n \in \N$ and $\sigma$ is a permutation on $\N$, is an isometry.
\end{lemm}

Let $|.|_D$ denote the equivalent Day's norm on $c_0$, that is
for $x=(x_n)_n \in c_0$,
$$|x|_D=\sup(\sum_{i=1}^k x_{n_i}^2/4^i)^{1/2},$$
where the sup is taken over $k \in \N$ and all $k$-tuples $(n_1,\ldots,n_k)$ of distincts elements of $\N$.
let $\norm{.}$ denote the corresponding norm on $X=c_0(Y)$, therefore for
$x=(y_k)_k \in X$,
$$\norm{x}=\sup(\sum_{i=1}^k \norm{y_{n_i}}^2/4^i)^{1/2},$$
and let $\norm{.}_{\infty}$ denote the sup norm on $X$,
$\norm{x}_{\infty}=\sup_k \norm{y_k}.$
Note that isomorphisms associated to a permutation on $\N$ and a sequence of signs are isometries on $X$ for $\norm{.}$. It remains to prove that $\norm{.}$ is LUR. Let $x=(y_k)_k \in X$ and $x_n=(y_{n,k})_k \in X$ be such that
$$\lim_n \norm{x_n}=\norm{x} \eqno (17)$$
and
$$\lim_n \norm{x+x_n}=2\norm{x}. \eqno (18)$$
We need to prove that $\lim_n \norm{x-x_n}=0$ or equivalently
$\lim_n \norm{x-x_n}_{\infty}=0$. Since $(x_n)_n$ is arbitrary satisfying (17) and (18) it is enough to prove that some subsequence of $(x_n)_n$ satisfies $\lim_n \norm{x-x_n}_{\infty}=0$.

Since, by elementary properties of $|.|_D$,
$$\norm{x+x_n} = |(\norm{y_k+y_{n,k}})_k|_D
\leq |(\norm{y_k}+\norm{y_{n,k}})_k|_D \leq \norm{x}+\norm{x_n},$$
we deduce from (17) and (18) that
$$\lim_n |(\norm{y_k}+\norm{y_{n,k}})_k|_D=2|(\norm{y_k})_k|_D. \eqno (19)$$
Since $|.|_D$ is LUR on $c_0$, \cite{DGZ} Theorem 7.3 p 69, we deduce from (17) and (19) that $$\lim_n |(\norm{y_k}-\norm{y_{n,k}})_k|_D=0,$$ therefore
$$\lim_n \max_k |\norm{y_{n,k}}-\norm{y_k}|=0. \eqno (20)$$
For any $n \in \N$, let $k_n \in \N$ be such that
$$\norm{x-x_n}_{\infty}=\norm{y_{k_n}-y_{n,k_n}}. \eqno (21)$$
Note that if $\lim_n k_n=+\infty$, then
$\norm{x-x_n}_{\infty} \leq 2\norm{y_{k_n}} +\max_k |\norm{y_{n,k}}-\norm{y_k}|$ converges to $0$. So passing to a subsequence we may assume that $(k_n)_n$ is constant equal to some $k_0 \in \N$. If $y_{k_0}=0$ then by (20), $\lim_n y_{n,k_0}=0$ and $\lim_n \norm{x-x_n}_{\infty}=\lim_n \norm{y_{k_0}-y_{n,k_0}}=0$. Therefore we may assume that $y_{k_0} \neq 0$.

Let $m \in \N$ be such that $m \geq |\{i \in \N: \norm{y_i} \geq \frac{1}{2}\norm{y_{k_0}}\}|$. Let $\beta=\frac{1}{2}\frac{\norm{y_{k_0}}}{2^m}$.
We prove that for $n$ large enough, 
$$\norm{y_{k_0}+y_{n,k_0}} \geq \beta. \eqno (22)$$
Indeed if (22) is contradicted then it is easy to see by the expression of $|.|_D$ that we may assume that for all $n$,
$$\norm{x+x_n}^2 \leq \sum_{i=1}^{+\infty} \frac{\norm{y_{k_i^n}+y_{n,k_i^n}}^2}{4^i}+\beta^2,$$
for some sequence $(k_i^n)_{i \geq 1}$ of distinct integers different from $k_0$. Let $\epsilon$ be positive. By (20) we deduce, for $n$ large enough,
$$\norm{x+x_n}^2 \leq (1+\epsilon)4\sum_{i=1}^{+\infty} \frac{\norm{y_{k_i^n}}^2}{4^i}+\beta^2,$$
So
$$\norm{x+x_n}^2 \leq (1+\epsilon)4\sum_{i=1}^{+\infty} \frac{\norm{y_{j_i}}^2}{4^i}+\beta^2,$$
where $(\norm{y_{j_i}})_{i \geq 1}$ is a non-increasing enumeration of
$\{\norm{y_k}, k \neq k_0\}$. Passing to the limit in $n$ and $\epsilon$, and using (18), we deduce
$$4\norm{x}^2 \leq 4\sum_{i=1}^{+\infty} \frac{\norm{y_{j_i}}^2}{4^i}+\beta^2 \leq 4\sum_{i=1}^{m} \frac{\norm{y_{j_i}}^2}{4^i}
+ \norm{y_{k_0}}^2 \sum_{i=m+1}^{+\infty}\frac{1}{4^i}+\beta^2,$$
therefore
$$4\norm{x}^2+\frac{\norm{y_{k_0}}^2}{4^m} \leq 
4(\sum_{i=1}^{m} \frac{\norm{y_{j_i}}^2}{4^i}
+\frac{\norm{y_{k_0}}^2}{4^{m+1}})+\frac{\norm{y_{k_0}}^2}{3.4^m}
+\beta^2 \leq
4\norm{x}^2+\frac{\norm{y_{k_0}}^2}{3.4^m}
+\beta^2.$$
We deduce that $\frac{2}{3.4^m}\norm{y_{k_0}}^2 \leq \beta^2$, a contradiction. Therefore (22) is proved.
Now
$$2\norm{x}^2+2\norm{x_n}^2-\norm{x+x_n}^2=
2\sum_{i=1}^{+\infty} \frac{\norm{y_{l_i}}^2}{4^i}
+2\sum_{i=1}^{+\infty} \frac{\norm{y_{n,l_i^n}}^2}{4^i}-
\sum_{i=1}^{+\infty} \frac{\norm{y_{n,m_i^n}+y_{m_i^n}}^2}{4^i},$$
where $(\norm{y_{l_i}})_i$, $(\norm{y_{n,l_i^n}})_i$,
and $(\norm{y_{n,m_i^n}+y_{m_i^n}})_i$ are non-increasing enumerations
of $(\norm{y_k})_k$, $(\norm{y_{n,k}})_k$, and $(\norm{y_k+y_{n,k}})_k$,
respectively.
Therefore
$$2\norm{x}^2+2\norm{x_n}^2-\norm{x+x_n}^2 \geq
2\sum_{i=1}^{+\infty} \frac{\norm{y_{m_i^n}}^2}{4^i}
+2\sum_{i=1}^{+\infty} \frac{\norm{y_{n,m_i^n}}^2}{4^i}-
\sum_{i=1}^{+\infty} \frac{\norm{y_{n,m_i^n}+y_{m_i^n}}^2}{4^i}.$$
Since by (17) and (18), $$\lim_n 2\norm{x}^2+2\norm{x_n}^2-\norm{x+x_n}^2=0,$$ we deduce by \cite{DGZ} Fact 2.3 p 45 that
$$\forall i \in \N, \lim_n 2\norm{y_{m_i^n}}^2+2\norm{y_{n,m_i^m}}^2-
\norm{y_{n,m_i^n}+y_{m_i^n}}^2=0. \eqno (23)$$
Let $K \in \N$ be such that for $k>K$, $\norm{y_k} \leq \frac{\beta}{4}$. By (20), we have for $n$ large enough and
$k>K$,
$$\norm{y_k+y_{n,k}}
 \leq 2\norm{y_k}+\frac{\beta}{4} \leq \frac{\beta}{2}.$$
By (22) we deduce that for $n$ large enough, $k_0 \in \{m_1^n,\ldots,m_K^n\}$. There exists $i$ such that $k_0=m_i^n$ for infinitely many $n$'s. Therefore from (23) we deduce, passing to a subsequence,
$$\lim_n 2\norm{y_{k_0}}^2+2\norm{y_{n,k_0}}^2-\norm{y_{k_0}+y_{n,k_0}}^2=0.$$
Since the norm $\norm{.}$ on $Y$ is LUR, this implies by \cite{DGZ}Proposition 1.2 p 42 that
$\lim_n y_{n,k_0}=y_{k_0}$. Finally
$$\lim_n \norm{x-x_n}_{\infty}=\lim_n \norm{y_{k_0}-y_{n,k_0}}=0.$$
\pff

\paragraph{Acknowledgements} We thank G. Godefroy and G. Lancien
 for their help concerning LUR-renorming theorems for Banach spaces, and
 G. Godefroy and V. Zizler for information about the articles \cite{J} and \cite{St}.

\

\

Valentin Ferenczi,

Institut de Math\'ematiques de Jussieu,

Universit\'e Pierre et Marie Curie - Paris 6,

Projet Analyse Fonctionnelle, Bo\^\i te 186,

4, place Jussieu, 75252 Paris Cedex 05,

France.

\

e-mail addresses: ferenczi@ccr.jussieu.fr, ferenczi@ime.usp.br.

\

\

El\'oi Medina Galego,

Departamento de Matem\'atica,

Instituto de Matem\'atica e Estat\' \i stica,

Universidade de S\~ao Paulo.

05311-970 S\~ao Paulo, SP,

Brasil.

\

E-mail: eloi@ime.usp.br.

\end{document}